\documentstyle[amscd,amssymb,verbatim]{amsart}

\theoremstyle{plain}
\newtheorem{Thm}{Theorem} 
\newtheorem{Lem}{Lemma}
\newtheorem{Prop}{Proposition}
\newtheorem{Def}{Definition}
\newtheorem{Cor}{Corollary}
\newtheorem{Rem}{Remark}

\renewcommand{\rm}{\normalshape}
\newcommand{\m}[2]{{\bold Mat}
\begin{Sb}
r_{{#1}}
\end{Sb}
({#2})}
\newcommand{\dd}{\partial_1 \times \dots\times \partial_n}
\newcommand{\dtr}{\partial_1 \times \partial_2 \times \partial_3}
\newcommand{\ptr}{p_1 \times p_2 \times p_3}
\newcommand{\alphatr}{\alpha_1 \times \alpha_2 \times \alpha_3}
\newcommand{\Ltr}{\Lambda \times \Lambda \times \Lambda}
\newcommand{\mtr}{{\Bbb TM'}_{0}^{3}}
\newcommand{\D}{{\Bbb \Delta}}
\newcommand{\N}{{\cal N}}

\newcommand{\A}{{\Bbb A}}
\newcommand{\mm}{{\Bbb M}^n _g}
\newcommand{\mn}{{\Bbb M}_{g,n}}
\newcommand{\mt}{{\Bbb T}^n _g}
\newcommand{\tm}{{\Bbb TM}^n _g}
\newcommand{\tmm}{{\Bbb TM'}_g ^n}
\newcommand{\k}{{\bold k}}
\newcommand{\g}{\Gamma}
\newcommand{\gl}[1]{L
\begin{Sb}
{\bold Mat}
\begin{Sb}
r_1
\end{Sb}(A_1)
\dots
{\bold Mat}
\begin{Sb}
r_n
\end{Sb}(A_n)
\end{Sb}
({#1})}

\begin{document}

%TOPMATTER
 
\title[A definition of graph homology and graph K-theory of algebras ]
{A definition of graph homology and graph K-theory of algebras }

\date{} 
\author{M. V. Movshev} 
\address{School of Mathematics \\ 
IAS,Olden Lane \\
Princeton,NJ 08540} 
\email{movshev@@ias.edu}
\thanks{This work was supported in part by A. Sloan Doctoral Fellowship and was completed in IAS during Spring semester of 1998}

\dedicatory{}

%End of topmatter 

\maketitle

\begin{abstract} 
We introduce and study elementary properties of graph homology of algebras. This new homology theory  shares many  features of cyclic and Hochschild homology. We also define a graph K-theory together with an analog of Chern character.
\end{abstract}

\section{Definition of fat graphs}
\begin{Def}\label{D:rtn}
Consider a graph $\g$. We call a
{\bf flag } a germ of an edge adjacent to a vertex.
On the enclosed picture the reader see a
graph with one vertex $a$, one edge
$x$ and two flags $\Lambda$ and       
$\Theta$.
\end{Def}

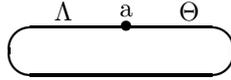
\setlength{\unitlength}{0.00083300in}
\begin{picture}(1530,464)(400,-957)
\thicklines
\put(3090,-794){\oval(1400,300)}
\put(3081,-676){\makebox(0,0)[lb]{$\bullet $ }}
\put(3081,-576){\makebox(0,0)[lb]{a }}
\put(3081,-1076){\makebox(0,0)[lb]{x }}
\put(2671,-601){\makebox(0,0)[lb]{$\Lambda$}}
\put(3451,-601){\makebox(0,0)[lb]{$\Theta$}}
\end{picture}

{\bf Convention} In these notes all graphs will be connected, non-directed with no vertices of valence one. When we say that a graph $\g$ is mapped to some topological space we mean a mapping of a CW-complex $CW(\g)$ associated to the graph.

\begin{Def}
A fat graph is a graph with a choice of a  cyclic order on the set of flags
adjacent to a vertex.
\end{Def}

   Let $X$ be a two dimensional oriented surface with $n \geq 1$ boundary components. Assume that a graph $\Gamma$, thought  as CW complex is embed into $X$ and is a deformation retract of the later. Then it inherits a fat structure from the surface: orientation of $X$ uniquely determine a cyclic order on flags at each vertex.
\begin{Prop} \label{P:wnu}
Let $D$ be a two or one dimensional disk. Let $A: D\rightarrow D$ be a homeomorphism which is the identity map on the boundary $\partial D$. Then the group of all such homeomorphisms is connected.
\end{Prop}
       It turns out. that one can recover the surface $X$ from the fat structure of $\g$:

  \begin{Prop} \label{P:vtr}
a) For every fat graph $\g$ there is an oriented surface with $n \geq 1$ boundary components, such that $i:CW(\g) \hookrightarrow X$ and $CW(\g)$ is a deformation retract of $X$. Moreover the fat structure of $\Gamma$ is induced from the orientation of $X$.
b) Let $i$ be  one such embedding. Any automorphism  of a fat graph can be extended to an orientation preserving homeomorphism of the  corresponding surface. The isotopy class of such  homeomorphism is fixed uniquely.
\end{Prop}
\begin{pf}
The proof is simple and is left to the reader as an exercise. In the proof of b) use proposition \ref{P:wnu}.
\end{pf}
\begin{Def}
Let $\g$ be a fat graph with an embedding of $i:CW(\g) \rightarrow X$ into an oriented surface with a boundary $X$. Suppose that $CW(\g)$ is a deformation retract of $X$ and a fat structure of $\g$ coincides with the one coming from the surface. We shall call a map $i$ a coherent embedding.  
\end{Def}

As a corollary of the last  proposition we get that the number of boundary components and the genus of the surface $X$ are invariants of the fat graph $\Gamma$. We denote them by $n(\Gamma)$ and $g(\Gamma)$ and call them the  number of boundary components and the genus of $\Gamma$.

We need to introduce some important groups. Fix a fat graph $\g$ and an oriented surface with a boundary $X$. Denote $\bar H(\g)$ the groups of homeomorphisms of $CW(\g)$.We have a subgroup $H(\g) \subset \bar H(\g)$ of homeomorphisms  which leave the fat structure invariant. Let $H_0(\g)$ be a connected component of the unit in $H(\g)$.
By aid of proposition \ref{P:wnu} one can easily see that there is a short exact sequence of groups:
\begin{equation}\label{E:sur}
1 \rightarrow H_0(\g) \rightarrow H(\g) \rightarrow \bar Aut(\g) \rightarrow 1
\end{equation}
The group $\bar Aut(\g)$ is the group of automorphisms of the fat graph $\g$ which preserve the fat structure.

Denote $\bar H(X)$ the group of all orientation preserving homeomorphisms of $X$. Let $H(X)$ be a subgroup of $\bar H(X)$  which consists of elements mapping each boundary component to itself. Denote by $H_0(X)$ the connected component of the unit in $\bar H(X)$. This data defines a pair of  short exact sequences 
\begin{equation}\label{E:fpr}
1 \rightarrow H_0(X)\rightarrow \bar H(X)\rightarrow \bar\g_{g}^n\rightarrow 1
\end{equation}
\begin{equation}\label{E:fyu}
1 \rightarrow H_0(X)\rightarrow H(X)\rightarrow \g_g^n\rightarrow 1
\end{equation}

\begin{Def}\label{D:yet}
We shall call $\bar \g_g^n$ a mapping class group and $ \g_g^n$ a mapping class group which preserves a marking.
\end{Def}

 Fix a coherent embedding $i: \g \rightarrow X$. Propositions \ref{P:wnu} and \ref {P:vtr} imply that there is a well defined homomorphism $h_i: \bar Aut(\g) \rightarrow \bar \g_g^n$.
\begin{Prop}
a) Homomorphism $h_i$ is injective. Denote $Aut(\g)=Im(h_i)\cap \g_g^n$. The definition of $Aut(\g)$ is independent on $i$. we shall call an element $a\in Aut(\g)$ a marking preserving automorphism.
\end{Prop} 
\begin{Def}\label {D:pyr}
we shall call an element $a\in Aut(\g)$ a marking preserving automorphism of a fat graph $\g$.
\end{Def}
\begin{Rem}\label{O:pot}
The short exact sequences \ref{E:fpr} and \ref{E:fyu} do not split. However the sequence \ref{E:sur} does. If we endow the edges of the graph $\g$ with a metric of a unit length. The group $Aut(\g)$  acts on $CW(\g)$ by isometries. Denote the isometry corresponding to $a\in Aut(\g)$ by $CW(a)$.
\end{Rem}

Two definitions are in order.
\begin{Def}\label{def:piu}
 We say that two maps $a,b :CW(\g) \rightarrow CW(\g')$ are isotopic if there are homeomorphisms $m \in H_0(\g)$ and  $n \in H_0(\g')$ such that the following diagram commutes:

$\begin{CD}
CW(\g)  @>\text{a}>>      CW(\g')\\
@VV\text{m}V   @VV\text{n}V \\
CW(\g)  @>\text{b}>>      CW(\g')
\end{CD} $
\end{Def}
\begin{Def}\label{D:tre}
 We say that two maps $i,i' :CW(\g) \rightarrow X$ are isotopic if there are homeomorphisms $m \in H_0(\g)$ and  $z \in H_0(X)$ such that the following diagram commutes:

$\begin{CD}
CW(\g)  @>\text{i}>>      X\\
@VV\text{m}V   @VV\text{z}V \\
CW(\g)  @>\text{i'}>>      X
\end{CD} $
\end{Def}

\subsection{Shrinking of edges} \label{S:tzp}
 Let $\g$ be a fat graph  with some coherent  embedding $i: CW(\g) \rightarrow X$ into an oriented surfaces with a boundary $X$ . Choose an edge $e$ of $\g$ which is not a loop. We shrink it, keeping the graph inside of the surface. As a result we get a new graph $\g/e$ with the fat structure induced from the ambient surface. By the construction we also get a coherent embedding $i'=e(i): CW(\g/e) \rightarrow X$which is defined uniquely up to isotopy. The number of edges and vertices of $\g/e$ is less by one than that of $\g$, but   $n(\g/e)=n(\g)$ and $g(\g/e)=g(\g)$. It is clear that the fat structure of $\g/e$ doesn't depend on the particular choice of $i$. The shrinking of an edge $e$ also defines   a CW map  $CW(e): CW(\g) \rightarrow CW(\g')$ defined up to isotopy. 

\begin{Prop}\label{E:pio}
Shrinking defines a one-to-one correspondence between isotopy  classes of coherent embeddings into $X$ of $\g$ and $\g/e$.
\end{Prop}
\begin{pf}
Exercise.
\end{pf}

\section{Definition of  category $\mm$} \label {S:wer}

The objects of $Ob(\mm)$ are the  fat graphs $\g$ such that $g(\g)=g$ and $n(\g)=n$ ($n\geq 1$, $g \geq 1$ or $n\geq 3$, $g=0$). All vertices of $\Gamma \in Ob(\mm)$ have a valence $\geq 3$.

 To  define the set  of  morphisms $Mor(\mm)$ we list generators and relations. Each automorphism $a \in Aut(\g)$ of a fat graph $\g \in Ob(\mm)$  which preserves a marking (see definition \ref{D:pyr} ) defines a generator in  $ Hom_{\mm} (\Gamma,\Gamma)$.  We denote it  by the same letter  $a$.

Similarly a  shrinking of an edge $e$ in a graph $\g$ defines a generator $e \in Hom(\g,\g/e)$.

In remark \ref {O:pot} and section \ref{S:tzp} we defined a  map $CW$ which assigns to any automorphism of a graph or to a  shrinking of an edge a mapping of appropriate $CW$ complexes. We use it to define the set of relations.

We say that two composed  morphisms $a,b$ between graphs $\g$ and $\g'$ are equal if the maps $CW(a)$ and $CW(b)$ are isotopic in a sense of definition \ref{def:piu}.

Introduce two sets $E(\g)$ and $V(\g)$ which are sets of edges and vertices respectively. A morphism $a$ in $\mm$, represented as a composition of automorphisms and shrinkings defines a map $a:E(\g) \cup V(\g) \rightarrow E(\g') \cup V(\g')$. Observe  that two morphisms $a,b$ are equal iff corresponding maps $a,b:E(\g) \cup V(\g) \rightarrow E(\g') \cup V(\g')$ are equal.

\section{The homotopy type of the nerve of  $\mm$}
\begin{Def}
Let $C$ be a small category. Denote $\N(C)$ the nerve of $C$ 
\end{Def}

we shall establish the following theorem:
\begin{Thm}\label{T:nsa}
Let $\g_g^n$ be a mapping class group of a surface of genus $g$ with $n$ boundary components which preserves a marking (see definition \ref{D:yet}). Let  $K(\g_g^n,1)$ be its classifying space.  Then there is a homotopy equivalence of topological spaces $\N (\mm)=K(\g_g^n,1)$
\end{Thm}
 
\subsection{Coverings  of the  nerve of a category}
Let us remind one general construction from \cite{Q}.
\begin{Def} \label{D:pre}
Let $C$ be a small category, $Set$ be a category of sets.  Let $f$ be a covariant functor from $C$ to  $Set$ which maps all morphisms of $C$ into isomorphisms. Define a new category $(C,f)$, whose objects are pairs $(A,x)$ ($A \in Ob(C)$ and $x \in f(A)$). A  morphism $m$ between $(A,x)$ and $(B,y)$ is  a morphism $\psi:A \rightarrow B$ such that $\psi(x)=y$. There is a functor  $L:(C,F)\rightarrow C$ , which forgets about the second factor.
\end{Def}
\begin{Prop} \cite{Q} \label {P:qlr}
$\N(L):\N(C,f) \rightarrow \N(C)$ is a covering map. The fiber is isomorphic to $f(A)$ for some $A \in Ob(C)$
\end{Prop}
{\bf Convention.} We say that two small categories $C$ and $C'$ are homotopic if their nerves are.Similarly two functors $f_1,f_2: C \rightarrow C'$ are homotopic if the corresponding maps of nerves are.

\subsection{Construction of the functor $f$} \label{S:pwn}
Fix an oriented surface $X$ of genus $g$ with $n$ boundary components .  For an object $\g \in \mm$ let us let $f(\g)$ to be a set of isotopy classes of coherent embeddings $i:\g \hookrightarrow X$.

It is clear that automorphisms of a fat graph $\g$ induce automorphisms of $f(\g)$. According to proposition \ref{E:pio} collapsing of an edge $e$ of $\g$ produces a one-to-one map $e:f(\g) \rightarrow f(\g/e)$. 
We  conclude that  the functor $f$ defines a covering of $\N (\mm)$. 
Denote he category $(\mm,f)$ by $\mt$.

{\bf Observation} Mapping class group $\g_g^n$ acts on $f$ by natural transformations.

\begin{Prop}\label{P:aie}
The action of $\g_g^n$ is free on $f(\g)$ for any $\g$
\end{Prop} 
\begin{pf}
Fix a fat graph $\g$

We check a  transitivity first. Each embedding $i:\g\hookrightarrow X$ ($i \in f(\g)$) defines some   cell decomposition $\Xi(i)$ of $X$. Here is it construction: each vertex $v \in \g$ of valence $n$ defines a two- cell $\xi_i(v)$ in a shape of an $n$-gon sitting inside of $X$ and intersecting the boundary $\partial X$ exactly in $n$ points. Each edge $e$ defines a cell in a shape of a rectangle $\xi_i(e)$, intersecting the boundary by a pair  opposite sides, the other pair sits on a boundary of some $n$-gons, corresponding to vertices which $e$ connects. One-cells are intersections of rectangle cell with a boundary and intersections of rectangle cells with n-gons. Zero-cell are defined similarly.

Suppose we have two embeddings $i,i'$ where $i,i' \in F(\g)$, which define two cell decompositions $\Xi(i)$ and $\Xi(i')$ of $X$. Define a map $z_{i.i'}$ which maps $\xi_i(v)$ into $\xi_{i'}(v)$ for each vertex $v$. Do similarly for edges. It is clear that the map $z_{i.i'}$  is defined correctly (it is only needed to check a consistency  on the boundary of the cells, but it is more or less obvious).

A triviality of a stabilizer.  Suppose  we have an embedding $i$, such that for some element $z\in \g_g^n$ the embedding $z(i)$ is isotopic to $i$. Without loss of generality we may assume that $\Xi(i)=\Xi(z(i))$ and $z$ is the identity on the one-skeleton. For each two-cell the map $z$ defines it homeomorphism, which is an identity on the boundary. It is known that the group of such homeomorphisms is connected, therefore the stabilizer of $\g_g^n$-action  is trivial.
\end{pf}

\begin{Lem} \label{L:pew}
For any pair of objects $(\g,i)$ and $(\g',i')$ the set\\ $Hom_{\mt}((\g,i),(\g',i'))$ is either empty or contains only one element.
\end{Lem} 
\begin{pf}
Let $a,b \in Hom_{\mm}(\g,\g')$ be two morphisms such that $f(a)(i)=f(b)(i)=i'$(see section \ref{S:pwn} for definition of $f$). If we choose  representatives in the isotopy class we get $f(a)(i)=i'$, $f(b)(i)=i''$ where $i'$ and $i''$ are isotopic (see definition \ref{D:tre}).

So we have a commutative diagram

$\begin{CD}
X  @>\text{L}>>      X\\
@AA\text{$i'\circ a$}A   @AA\text{$i''\circ b$}A \\
CW(\g)  @>\text{l}>>      CW(\g')
\end{CD} $

In this diagram $L \in H_0(X)$, $l \in  H_0(\g)$. Since $i'$ and $i''$ are embeddings there is an element $l' \in \bar H(\g')$, which makes the following diagram commute:

$\begin{CD}\label{D:ope}
X  @>\text{L}>>      X\\
@AA\text{i'}A   @AA\text{i''}A \\
CW(\g')  @>\text{l'}>>      CW(\g')\\
@AA\text{ a}A   @AA\text{ b}A \\
CW(\g)  @>\text{l}>>      CW(\g)
\end{CD} $

The homeomorphism $l'$ leaves the fat structure of $\g'$ invariant and must belong to $H(\g')$. Moreover $l' \in H_0(\g)$ otherwise it would contradict to commutativity of the top part of the last  diagram.Commutativity of the low part of that diagram implies that by definition of relations in the category given in section \ref{S:wer}  $a=b$ as a morphism in $\mm$.
\end{pf}

\subsection{Definition of a category $\A$}
We introduce a category $\A$ which will be used in the proof of theorem \ref{T:nsa}.
Denote $\bar X$ to be a compact surface such that $X \subset \bar X$ and $\bar X \backslash X$ is a union of discs. 

\begin{Def}\label{D:pbe}
  An object $\Sigma \in Ob(\A)$ will be a cell decompositions of $\bar X \supset X$. Any such cell decomposition will have as  one-skeleton  the image of $i:\g \rightarrow X \subset \bar X$ for some fat graph $\g$ and coherent embedding $i$. The union of  two-cells is  the complement to  the image of $i$. Two cell decompositions $\Sigma$ and $\Sigma'$  will be considered as equal if they are isotopic.

 The morphisms in $\A$ are shrinkings of one-cells which are not loops.
\end{Def}

The nerve of this category has been studied in \cite{H}.
To make a link with his notation we need the following 

\begin{Lem}
Suppose a graph $\g \subset \bar X$  is a one-skeleton of a cell decomposition of $\bar X$(that is the complement of the graph is a union of disks).  Then the dual graph $\g\spcheck$ also defines a cell decomposition of $\bar X$.
\end{Lem}

\begin{pf} Let $i:\Lambda \hookrightarrow \bar X$ an embedding of some graph $\Lambda$ into  $\bar X$.
Let  $v(\Lambda)$ be the number of vertices in $\Lambda$,  $e(\Lambda)$ be the number of edges in $\Lambda$, and $c(\Lambda)$ be the number of connected components in the complement $\bigcup_{i=1}^{c(\Lambda)} X(\Lambda)_i=\bar X \backslash \Lambda$. So  in our notations $\bar X= \Lambda \cup  \bigcup_{i=1}^{c(\Lambda)} X(\Lambda)_i$. By aid of Mayer-Vietoris long exact sequence we get an equation for Euler characteristics:
\begin{equation}
\chi(\bar X)=\chi(\Lambda)+\bigoplus \chi(X_i(\Lambda))=v(\Lambda)-e(\Lambda)+\bigoplus \chi(X_i(\g))
\end{equation}
  
For graph $\g$ and its dual $\g\spcheck$ we have the following identities: 
\begin{equation}
v(\g)=c(\g\spcheck),\quad e(\g)=e(\g\spcheck), \quad v(\g\spcheck)=c(\g)
\end{equation}

  Now let us use the fact that $X_i(\g)$ are contractible . It is equivalent to  $\chi(X_i(\g))=1$.  Then we have 
\begin{align}
&\chi(X)=v(\g)-e(\g)+\bigoplus_{i=1}^{c(\g)} \chi(X_i(\g))= \notag\\
=&v(\g\spcheck)-e(\g\spcheck)+\bigoplus_{i=1}^{c(\g\spcheck)}\chi(X_i(\g\spcheck))=\notag\\
=&\bigoplus_{i=1}^{c(\g)}  \chi(X_i(\g))-e(\g)+\bigoplus_{i=1}^{v(\g)} \chi (X_i(\g\spcheck)) \notag
\end{align}

 It implies that $\bigoplus_{i=1}^{ v(\g)}\chi(X_i(\g\spcheck))=v(\g)$ which in its turn gives $\chi(X_i(\g\spcheck))=1$.
\end{pf}

Fix $n$ points on a compact surface $X$. Define a rank k- arc-system as a collection of $k+1$ paths on the surface $X$ connecting points, such that they intersects only at the endpoints and no two of them are isotopic. We say that the arc system fill the surface iff the complement to the set of paths is a collection of disks.

Observe that a collapsing of a one-cell in the language of a cell decomposition $\Sigma$ of  $\bar X$ in the dual language corresponds to removing of a path in the arc-system made out of the $\g\spcheck$ ($\g\spcheck$ is a one-skeleton of a dual cell decomposition of X). We have a simplicial set, we call it $A$ , whose $k$ simplices are rank-k ark-systems and the faces maps are removals of arcs. Denote $A_{\infty}$ the simplicial set of arc-systems which do not fill the surface. Let $Y$ be a simplicial subset of the first barycentric subdivision $A^0$ of $A$ which has no faces in $A_{\infty}$.

 \begin{Thm} \cite {H}\label {P:rew}
The geometric realization $\Delta(Y)$ of $Y$ is contractible .
\end{Thm}

\begin{Prop} \label {P:vwt}
  The nerve $\N(\A)$ is identically  equal to $\Delta(Y)$. 
\end{Prop}
\begin{pf}
Exercise.
\end{pf}

\subsection{Sufficient conditions for a functor to  induce a homotopy equivalence of nerves}

Suppose we are given two small categories $C$ and $C'$ and a functors $g_1,g_2:C \rightarrow C'$.

We formulate a condition when two such functors are homotopic.
\begin{Prop}\cite{Q} \label{P:kgd}
Suppose there is a morphism between functors $g_1$ and $g_2$. Then they are homotopic. 
\end{Prop}
Fix a functor $g:C \rightarrow C'$.
\begin{Def}
For any object $Y \in Ob(C')$ a category $Y\backslash g$ is defined  as follows: its objects are  diagrams $Y \overset{v}{\rightarrow} g(X)$ for some object $X\in C$. Morphism are $w:X \rightarrow X'$ such that $v'=g(w)v$.
\end{Def}
\begin{Prop} \cite{Q}\label{P:car}
If $Y\backslash g$ is contractible for any $Y \in Ob(C')$ then $g$ is a homotopy equivalence.
\end{Prop}

\begin{Def}
There is also a dual construction. For any object $Y \in Ob(C')$ a category $ g/Y$ as follows: its objects are  diagrams $ g(X) \overset{v}{\rightarrow}Y $ for some object $X\in C$. Morphisms are $w:X \rightarrow X'$ such that $v=v'g(w)$
\end{Def}

\begin{Prop} \cite{Q} \label{P:mar}
If $g/Y$ is contractible for any $Y \in Ob(C')$ then $g$ is a homotopy equivalence
\end{Prop}
\begin{Prop} \cite{Q} \label{P:bgt}
A category $C$ having initial or terminal object is contractible.
\end{Prop}
\subsection{Proof of  theorem \ref {T:nsa}} 
\begin{pf}

we shall show that the  universal cover of  $\N(\mm)$ is contractible and $\pi_1(\N(\mm))=\g_g^n$. we shall identify the universal cover of $\N(\mm)$ with $\N(\mt)$. According to proposition \ref{P:qlr} we have a covering map $\N(L): \N(\mt) \rightarrow \N(\mm)$. All we need is to check  that $\N(\mt)$ is contractible.

There is a forgetting functor $g:\mt \rightarrow \A$, which assigns to a pair $(\g,i)$ the cell decomposition of $\bar X \supset X$  described in the definition \ref{D:pbe}. According to theorem \ref{P:rew} and  proposition \ref{P:vwt} the category $\A$ is contractible. To prove that $\mt$ is contractible we show that $g$ is a homotopy equivalence.

We employ proposition \ref {P:car} and show that for any object  $\Sigma \in \A$ the category $\Sigma\backslash g$ is contractible .

 Fix an object $\Sigma \in Ob(\A)$. Suppose we have an object in  $Ob(\Sigma\backslash g)$ which is represented by a diagram:
\begin{equation}
 \Sigma \overset{\psi}{\rightarrow} g((\g',i'))
\end{equation}
In more condensed notations we shall denote this diagram as $\{\psi,(g',i')\}$
 We can choose an object $(\g,i) \in Ob(\mt)$ and a morphism $\bar \psi \in Hom_{\mt}((\g,i),(\g',i'))$ such that $g(\g,i)=\Sigma$ and $g( \bar \psi)=\psi$.

 The object $(\g,i)$ is some coherent embedding of $\g$ into $X$ which gives rise to $\Sigma$. The set of such embeddings is finite and each element has a form of composition $i\circ a$ , where $a$ is an automorphism of $\g$. 

We fix one of such $(\g,i)$. The  arguments given above imply that for any   $\{\psi,(\g',i')\}$ there is a morphism $b:\{id,(g,i)\} \rightarrow \{\psi,(\g',i')\}$. By lemma \ref{L:pew} this isomorphism is unique. We can claim that the category  $\Sigma\backslash g$ has an initial object. According to proposition \ref{P:bgt} the nerve  $\N(\Sigma\backslash g)$  is contractible.

In proposition \ref{P:aie} we established that the group $\g_g^n$ acts by natural transformations of the functor $f$ and the action is free on each $f(\g)$. This implies that the action commutes with the projection $\N(L):\N(\mt) \rightarrow \N(\mm)$ and is free on the fibers. We can conclude that $\N(\mm)$is indeed $K(\g_g^n,1)$ 
\end{pf}

\section{Some modifications of the category $\mm$}

We introduce a new category $\tmm$.The fat graphs which belong to $Ob(\tmm)$ are allowed to have two-valent vertices.  As in case of $\mm$ we define morphisms- generators  first. These will be shrinkings of edges(which are not loops), automorphisms (which preserve the fat structure and marking) and insertions of points on edges. The set of relations defined as follows.
Two morphisms $\alpha,\beta:\g \rightarrow \g'$ are equal iff $CW(\alpha)$ and $CW(\beta)$ are homotopic relative to zero-skeleton of $CW(\g)$.

The category $\tm$ is a modification of $\tmm$.  In the definition of morphisms  we omit  insertions of vertices 
\begin{Rem}
In the definition of $\tm$ $\tmm$  we can relax a bit restrictions on $n,g$ and allow all $n \geq 1$ $g \geq 0$.
\end{Rem} 
\begin{Prop}\label{P:wbi}
a) There is an inclusion functor  $i:\mm \rightarrow \tmm$, whose image is a full subcategory of $\tmm$.

b)The functor $i$ is a homotopy equivalence.
\end{Prop}
\begin{pf}
We leave the proof of a) to the reader.

b) we shall construct a functor $er:\tmm \rightarrow \mm$ such that $\N(er)$ is a  homotopy inverse  to $\N(i)$.  The functor $er$ erases all vertices of valence two. The composition $er\circ i$ is isomorphic to identity functor of $\mm$. There is a morphism $i\circ er \rightarrow id$, which inserts all erased vertices. By proposition \ref{P:kgd}  it  proves the claim.  
\end{pf}

\subsection{A homotopy equivalence of $\tmm$ and $\tm$}

We start with some preliminaries.
\begin{Prop}\label{P:qmz}\cite {Q}
Suppose we have a direct system of categories $f_i:C_i \rightarrow C_{i+1}$ such that $f_i$ are  homotopy equivalences. The category $C$ is a direct limit  $C=\underset{\rightarrow}{\lim}C_i$. Then the functor $C_i \rightarrow C$ is a homotopy equivalence.
\end{Prop} 

\begin{Def}
Denote $[n]$ an ordered set $1<2<\dots <n$.
A category $\D$ has its objects $\{[n]| n \geq 1\}$. The morphisms are surjective nondecreasing maps.
\end{Def}
\begin{Def}
Denote by $\D_2$ a full  subcategory of $\D$ whose objects are $\{[n]| n \geq 2\}$.
\end{Def}
\begin{Prop} \label{P:uew}
Fix an object $[n] \in Ob(\D_2)$. Then the category $(id/[n])$ is contractible.
\end{Prop}
\begin{pf}
Indeed the category $(id/[n])$ has  a terminal object $[n] \overset{id}{\rightarrow}[n]$.
\end{pf}

\begin{Prop} \label{P:dnq}
Consider a full subcategory $\D_2(n)=\{[k]|2 \leq k \leq n\} \subset \D_2$. This category is contractible for all $n \geq 2$. 
\end{Prop}
\begin{pf}
The category $\D_2(n)$ sits inside $\D(n)=\{[k]|1 \leq k \leq n\}$. The nerve of $\D(n)$ is the first barycentric subdivision of the standard $n$-simplex $\bold \D$. The nerve of $\D_2(n)$ sits inside of  $\bold \D$ and is the union of all closed cells which have no common points with the $0$-skeleton of  $\bold \D$(before subdivision). If we place the standard $n$-simplex into ${\Bbb R}^n$ so that the mass center is at the origin then a homotety with the  coefficient $t$ ($0\leq t\leq 1$)  preserves $\N(\D_2(n))$. 
\end{pf}
\begin{Prop}\label{P:mbs}
The category $\D_2$ is contractible.
\end{Prop}
\begin{pf}
Apply propositions \ref{P:qmz} and  \ref{P:dnq}.
\end{pf}

\begin{Thm}\label{T:poz}
The inclusion functor $i:\tm \rightarrow \tmm$ satisfy the following property: for every $\g \in \tmm$ the category $(\g\backslash i)$ is contractible. Due to proposition \ref {P:car} $i$ is a homotopy equivalence.
\end{Thm}

\begin{pf}
It is easy to see that every morphism $\phi \in Hom_{\tmm}(\g,\g')$ can be {\it  uniquely}  factored $\g \overset{\phi'}{\rightarrow} \Phi \overset{\alpha}{\rightarrow} \g'$ such that  $\phi' \in Hom_{\tm}(\g,\Phi) \subset  Hom_{\tmm}(\g,\Phi)$ and $\alpha$ makes only insertions of vertices.
\begin{Lem}
The correspondence $(\g \overset{\phi}{\rightarrow}\g') \leadsto (\g \overset{\phi}{\rightarrow}\Phi) $ is actually a functor , which we denote by $er:(\g\backslash i) \rightarrow (\g \backslash id_{\tm})$.
\end{Lem}
\begin{pf}
The proof readily follows from uniqueness of factorization.
\end{pf}

The category $(\g \backslash id_{\tm})$ is contractible, since it has an initial object $\g  \overset{id}{\rightarrow} \g$. All we need is to show that $er$ is a homotopy equivalence.
\begin{Def}
Let $f:C \rightarrow C'$ be a functor between two small categories. Fix $Y \in Ob(C')$. Define a category $f^{-1}(Y)$ which consists of all $X$ such that $f(X)=Y$, the morphism are those maps which cover the identity $id \in Hom_{C'}(Y,Y)$. There are inclusions of categories $f^{-1}(Y) \rightarrow (Y\backslash f)$ ($X \leadsto Y\overset{id}{\rightarrow} f(X) $) and  $f^{-1}(Y) \rightarrow (f/Y)$ ($X \leadsto f(X)\overset{id}{\rightarrow} Y$).
\end{Def}

Our plan is to show $a)$ that for every $Y \in Ob(\g \backslash id_{\tm})$  the category $er^{-1}(Y)$ is contractible and $b)$the functor $g: er^{-1}(Y) \rightarrow (Y \backslash er)$ is a homotopy equivalence. These facts in combination with proposition \ref {P:car} yield the proof.

$a)$ Let us set the notations. The object  $Y \in Ob(\g \backslash id_{\tm})$ is a diagram $\g \overset{\phi}{\rightarrow}\Phi$ ($\phi \in Hom_{\tm}(\g,\Phi)$).  It is clear that every object of  $er^{-1}(Y)$ is represented by a diagram $\g\overset{\phi}{\rightarrow}\Phi \overset{\alpha}{\rightarrow}\g'$. The morphism $\alpha$ inserts vertices on edges of $\Phi$. Let $E(\Phi)$ be a set of edges of the  graph $\Phi$. To each element $e \in E(\Phi)$ we assign a weight  by the rule : $w_{alpha}(e)=$ the number of vertices inserted on $e$ by $\alpha$ plus two. Suppose $\#(E(\Phi))=s$. We assign 

\begin{equation}
\g \overset{\phi}{\rightarrow}\Phi \overset{\alpha}{\rightarrow}\g'  \leadsto [w_{\alpha}(e_1)] \times \dots \times [w_{\alpha}(e_1)] \in Ob(\D_2^{\times s})
\end{equation}

It is easy to see that this is an isomorphism of categories. By proposition \ref {P:mbs} the category $er^{-1}(Y)$ is contractible.

$b)$ It is easy to see that every object $X$ of the category $(Y\backslash er)$  is represented by a commutative  diagram 

\begin{equation}\label{D:ywk}
\begin{CD}
\Psi'\\
@AA{\alpha'}A\\
\Phi' @<\gamma<<\Phi\\
@A{\phi'}AA    @AA{\phi}A\\
\g @= \g\\
\end{CD}
\end{equation}

Here $\g,\phi, \Phi$ are fixed, $\gamma,\phi' \in Mor(\tm)$, $\alpha$ is insertion of vertices. we shall show that for every $X \in Ob(Y \backslash er)$ the category $(g/X)$ is contractible. With proposition \ref {P:mar} it enough for the  proof that $g$ is a homotopy equivalence.

We identify objects of  this  category with completions of upper right corner  of the diagram  \ref{D:ywk} :
\begin{equation}\label{D:wpb}
\begin{CD}
\Psi'@<\delta<< \Psi\\
@A{\alpha'}AA  @AA{\alpha}A\\
\Phi' @>\gamma>>\Phi\\
@A{\phi'}AA    @AA{\phi}A\\
\g @= \g\\
\end{CD}
\end{equation}  

``Free variables `` here are morphisms $\alpha$ (insertions of vertices) and $\delta$ $\in$\\ $ Hom_{\tm}(\Psi,\Psi')$. 

To establish that this category coincides with the one which is obviously contractible we elaborate one the structure of the set of objects.
The morphism $\gamma$ contracts some edges of $\Phi$ so we have an inclusion $E(\Phi') \subset E(\Phi)$. Suppose the set $E(\Phi)$ has $p$ elements, so $e_1,\dots,e_p \in E(\Phi')$, $e_{p+1},\dots,e_s \in E(\Phi)\backslash E(\Phi')$.
Introduce a category $(id_{\D_2}/[w_{\alpha'}(e_1)]) \times \dots \times (id_{\D_2}/[w_{\alpha'}(e_p)])\times \D_2^{\times(s-p)}$. As in the prove of statement $a)$ it is easy to see that $(g/X)$ is isomorphic to this category. But the later is contractible due to propositions \ref {P:uew} and \ref{P:mbs}.
\end{pf}
\subsection{Homology of categories}\label{S:nnb}

Let $f:C \rightarrow C'$ be a functor between two small categories and let $G:C' \rightarrow Vect$ be a covariant functor from $C'$ to the category of linear spaces. 

There is a map
\begin{equation}\label{E:pod}
\underset{\overset{\rightarrow}{C}}{\lim}G \circ f \rightarrow \underset{\overset{\rightarrow}{C'}}{\lim}G 
\end{equation}
Denote higher derived functors $R^i\underset{\overset{\rightarrow}{C'}}{\lim}G 
\overset{def}{=}H_i(C',G)$ -the homology of $C'$ with coefficients in $G$. The map \ref{E:pod} can be extended to a map of higher homology:
\begin{equation}\label{E:pld}
H_i(C,G \circ f) \rightarrow H_i(C',G)
\end{equation}

\subsection{Sufficient condition for a functor to induce homology equivalence}\label{S:neu}
We are interested in sufficient condition when maps \ref{E:pld} are isomorphisms. 

\begin{Def}
We say that a functor $f$ is a homology equivalence between if the map \ref{E:pld} is an isomorphism for any functor $G$ 
\end{Def}

We describe $E^1$ term of  a general spectral sequence which measure a deviation of the map  \ref{E:pld} from being isomorphism. The spectral sequence converges to $H_.(C,G \circ f)$.

  The group $E_{i,j}^1$ is equal 
\begin{equation}
E_{i,j}^1=\bigoplus_{[X_1\leftarrow \dots \leftarrow X_i]} [X_1\leftarrow \dots \leftarrow X_i]\otimes H_j((X_1\backslash f))\otimes G(X_i)
\end{equation}

The differential is equal $d=\sum_1^{i-1} (-1)^{k}d_k$. Maps $d_k$ ($k \neq 1,i-1$) replace a symbols $[X_1\leftarrow \dots \leftarrow X_{k-1}\overset{\alpha}{\leftarrow} X_k\overset{\beta}{\leftarrow} X_{k+1}\leftarrow \dots \leftarrow X_i]$ by $[X_1\leftarrow \dots \leftarrow X_{k-1}\overset{\alpha \circ \beta}{\leftarrow} X_{k+1}\leftarrow \dots \leftarrow X_i]$.

The map 

\begin{align}
d_{i-1}:&[X_1\leftarrow \dots \leftarrow X_{i-1} \overset{\alpha}{\leftarrow} X_i]\otimes H_j((X_1\backslash f))\otimes G(X_i)\\
\overset{id \otimes G(\alpha)}{\rightarrow}& [X_1\leftarrow \dots \leftarrow X_{i-1}]\otimes H_j((X_1\backslash f))\otimes G(X_{i-1})
\end{align}

The map 
\begin{align}
d_1:&[X_1\overset{\beta}{\leftarrow} X_2 \dots \leftarrow X_i]\otimes H_j((X_1\backslash f))\otimes G(X_i)\\
 \overset{H(\beta) \otimes id}{\rightarrow}& [X_2\leftarrow \dots \leftarrow X_{i}]\otimes H_j((X_2\backslash f))\otimes G(X_{i})
\end{align}

\begin{Rem}\label{R:ngc}
It is obvious that if homology $ H_.((X\backslash f))$ do not depend on the object $X$ and all maps $H(\beta)$ are isomorphisms  then we can push our calculations a little further and conclude that $E^2_{.,.}$ is equal to $ H_.((X\backslash f))\otimes H_.(C',G)$.
\end{Rem}

Suppose that $ H_i((X\backslash f))=0$ for $i>1$. Then the spectral sequence degenerates and  we proved the following result:
\begin{Prop} \label{P:rte}
Suppose that  for every $Y \in Ob(C')$ the category  $(Y\backslash f)$ is contractible. Then the maps \ref{E:pld} are isomorphisms.
\end{Prop}

We want to apply this sufficient condition to the inclusion functor $i:\tm \rightarrow \tmm$
\begin{Prop}\label{P:ghe}
For any functor $G:\tmm \rightarrow Vect$ the groups $H_i(\tm,G \circ i)$ and  $H_i(\tmm,G)$ are isomorphic.
\end{Prop}
\begin{pf}
In theorem \ref{T:poz} we showed that the categories $(\g \backslash i)$ are contractible. The claim follows from proposition \ref{P:rte}.
\end{pf}

\section{Boundary functors}
 I owe to M. Kontsevich the idea of these  functors.

\subsection{Functor of flags etc.}
Fix a fat graph $\g$ . Let $fl(\g)$ be a set of flags as it defined in \ref {D:rtn}. Introduce a transformation $D$ acting on $fl(\g)$ (c.f. \cite{K}). Fix a flag $f_1$ adjacent to a vertex $v$ .Let $f_2$ be a flag which stays right ahead of $f_1$ in the cyclic order of $v$. Suppose the flag $f_2$ is one end of the edge $e_{v,w}$. The other one is $f_3=D(f_1)$. We tried to illustrate  what are the geometric images staying behind this construction on the figure \ref{P:lkj}.
 Any fat graph can be coherently immersed into ${\Bbb R}^2$. The later is equipped with the standard  orientation. It is easy to see that this immersion is unique up to isotopy, if isotopy is allowed to bring edges over vertices.
Fattening the edges we locally get a surface which may look like this:

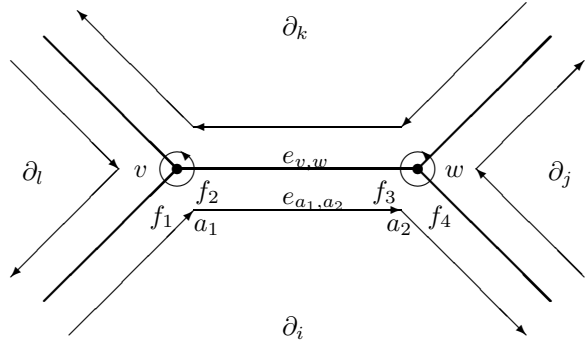
\begin{equation}\label{P:lkj}
\end{equation}
\setlength{\unitlength}{0.00057300in}
\begin{picture}(5199,3024)(0,-3298) 
\thicklines
\put(2701,-1786){\line( 1, 0){2175}}
\put(4876,-1786){\line( 1, 1){1200}}
\put(4876,-1786){\line( 1,-1){1200}}
\put(2701,-1786){\circle*{100}}
%\put(3751,-1786){\circle*{100}}
\put(4876,-1786){\circle*{100}}
\put(1501,-2986){\line( 1, 1){1200}}
\put(1501,-586){\line( 1,-1){1200}}
\thinlines
\put(1726,-3286){\vector( 1, 1){1125}}
\put(5401,-1786){\line( 1,-1){975}}
\put(5401,-1786){\vector(-1, 1){0}}
\put(1801,-361){\line( 1,-1){1050}}
\put(1801,-361){\vector(-1,1){0}}
\put(2851,-2161){\vector( 1, 0){1875}}
\put(4726,-2161){\vector( 1,-1){1125}}
\put(2851,-1411){\line( 1, 0){1875}}
\put(2851,-1411){\vector(-1, 0){0}}
\put(4726,-1411){\line( 1, 1){1125}}
\put(4726,-1411){\vector(-1,-1){0}}
\put(5401,-1786){\vector( 1, 1){975}}
\put(1201,-2761){\line( 1, 1){975}}
\put(1201,-2761){\vector(-1,-1){0}}
\put(1201,-811){\vector( 1,-1){975}}
\put(3651,-1736){\makebox(0,0)[lb]{$e_{v,w}$}}
\put(3651,-2100){\makebox(0,0)[lb]{$e_{a_1,a_2}$}}
\put(2730,-1620){\vector(-1, 1){0}}
\put(2701,-1786){\circle{300}}
\put(4905,-1620){\vector(-1, 1){0}}
\put(4876,-1786){\circle{300}}
\put(2310,-1860){\makebox(0,0)[lb]{$v$}}
\put(5120,-1860){\makebox(0,0)[lb]{$w$}}
\put(2851,-2350){\makebox(0,0)[lb]{$a_1$}}
\put(4590,-2350){\makebox(0,0)[lb]{$a_2$}}
\put(2451,-2281){\makebox(0,0)[lb]{$f_1$}}
\put(2861,-2061){\makebox(0,0)[lb]{$f_2$}}
\put(4460,-2061){\makebox(0,0)[lb]{$f_3$}}
\put(4960,-2281){\makebox(0,0)[lb]{$f_4$}}
\put(3651,-3286){\makebox(0,0)[lb]{$\partial_i$}}
\put(6051,-1886){\makebox(0,0)[lb]{$\partial_j$}}
\put(3651,-580){\makebox(0,0)[lb]{$\partial_k$}}
\put(1300,-1886){\makebox(0,0)[lb]{$\partial_l$}}
\end{picture}

On this picture thick lines correspond to a fragment of a  fat graph $\g$. It is clear that operator $D$ drags a flag along some boundary component. On the picture the boundary component corresponding to the flag $f_2$ is $\partial_i$.

  The set $fl(\g)$ brakes into the union of  orbits of $D$ action. Denote the set of orbits by $b(\g)$. 
\begin{Rem} \label{R:pbq}
The geometric meaning of the set $b(\g)$ is transparent. An orbit $a \in b(\g)$ defines an oriented path in $CW(\g)$. If one choose a coherent embedding $i:CW(\g) \rightarrow X$, then the path corresponding to the orbit  will be a image under deformation retraction of the one of the boundary components . The orientation of the boundary is induced by orientation of $X$. We conclude that for graph $\g$ with $n$ boundary components $\#\{b(\g)\}=n$ 
\end{Rem}

\begin{Rem}
Fix some automorphism $\alpha$ of a fat graph $\g$. The action of   $\alpha$ on $fl(\g)$ commutes with $D$, inducing the action on $b(\g)$. The automorphism $ \alpha$ preserves marking in a sense of definition \ref {D:pyr} iff its action on $b(\g)$ is trivial.
Let $\beta: \g \rightarrow \g'$ be an  insertion of a vertex or a shrinking of an edge. Though $\beta$ doesn't define a map of $fl(\g) \rightarrow fl(\g')$ it does define  $b(\g) \rightarrow b(\g')$. In other words we have a functor $b: \tmm \rightarrow Set$.
\end{Rem} 
The functor $b$ transforms morphisms into isomorphisms, therefore according to proposition \ref{P:qlr} defines a covering of the nerve $\N(\tmm)$. We know that the fundamental group of $\N(\tmm)$ is the mapping class group $\g_g^n$.
\begin{Prop}\label{P:tal}
The monodromy action of $\g_g^n$ on $b(\g)$ is trivial.
\end{Prop}

\begin{pf}
we shall examine a covering over $\N(\mm)$ defined by $b$ . Denote the category $(\mm, b)$ by $B$ ( see definition \ref {D:pre}). 

As we know there is a universal covering $\N(\mt)\rightarrow \N(\mm)$ with a fiber $f(\g)$-the set of all coherent embeddings of $\g$ into oriented surface $X$. This covering is a universal principal bundle and every covering is induced from it. On just have to fix some set with $\g_g^n$ action and carry out the induction. The suitable set is the set of boundary components with their natural orientation. The action of the mapping class group on this set is trivial. By definition  of the induced covering and remark \ref{R:pbq} this action coincides with the monodromy action.  
\end{pf}

\begin{Def}
A boundary component of the category $\tmm$ is a choice of an element  $o \in \underset{\overset{\longleftarrow}{\tmm}}{\lim} b$.
\end{Def}

 Proposition \ref{P:tal} guarantees that the set $\underset{\overset{\longleftarrow}{\tmm}}{\lim} b$  has exactly $n$ elements.

One can define the same way boundary components for categories $\mm$ and $\tm$.

\subsection{Boundary functors} 
    Our plan is to associate to an element $o_i \in \underset{\overset{\longleftarrow}{\tmm}}{\lim} b$  a functor  $ \tmm \rightarrow \Lambda$, which we denote $\partial_i$.
\begin{Def}\label{D:,qe}
 An element  $o_i$ defines an orbit $o_i(\g)$ in the set $fl(\g)$ of flags. Vertices of the polygon $\partial_i(\g)$ will be ordered triples $(f_1, f_2, v)$, where $f_2$ is a flag which belongs to the orbit $o_i(\g)$, $v$ is the vertex to which $f_2$ is adjacent, $f_1$ a flag which is behind of  $f_2$ in the cyclic order .  $a_1=(f_1, f_2, v)$ and  $a_2=(f_3, f_4, w)$ are connected by directed edge $e_{a_1,a_2}:a_1 \rightarrow a_2$ iff $f_2$ and $f_3$ are opposite ends of an edge $e_{v,w}$ connecting vertices $u$ and $v$. It is transparently clear $\partial_i(\g)$ is a {\it cyclic} graph. It is worthwhile to look  at the  picture \ref{P:lkj} and see how our formal definition identifies with simple geometric construction. 
\end{Def}
\begin{Prop}
The correspondence $\partial_i$ is a functor from $\tmm$ to $\Lambda$. The category $\Lambda={\Bbb TM'}_{0}^{2}={\Bbb TM'}_{0,1}$ is the cyclic category introduced by Connes (see \cite{C}).
\end{Prop}
\begin{pf}
Exercise.
\end{pf}

The functors $\partial_i(\g)$ for various $i$ will be the main source of graph objects, that is functors $f:\tmm \rightarrow Vect$. $Vect$ is  category of vector spaces. Indeed for each $i$ fix a functor $f_i:\Lambda \rightarrow Vect$. The tensor product 
\begin{equation}
f=\bigotimes f_i \circ\partial_i
\end{equation}
gives an example of a graph object.

\section{ A spectral sequence for graph homology}
In this section we elaborate on the spectral sequence introduced in the section \ref{S:neu} applied to a functor $f=\dd:\tmm \rightarrow \Lambda \times \dots \times \Lambda$.

Fix an object $[k_1]\times  \dots \times [k_n] \in Ob(\Lambda \times \dots \times \Lambda)$ . Our present aim is to understand the homotopy type  of the category $([k_1]\times  \dots \times [k_n]\backslash \partial_1 \times \dots \times \partial_n)$. 

Let $\phi:[k_1]\times  \dots \times [k_n] \rightarrow [k'_1]\times  \dots \times [k'_n]$ be a morphism. It defines a covariant functor 
\begin{equation}
\psi^*: ([k'_1]\times  \dots \times [k'_n]\backslash \partial_1 \times \dots\times  \partial_n) \rightarrow ([k_1]\times  \dots \times [k_n]\backslash \partial_1 \times \dots  \times\partial_n)
\end{equation}

\begin{Thm}\label{T:baw}
For any morphism $\psi$ the functor $\psi^*$ establish a homotopy equivalence between  $([k'_1]\times  \dots \times [k'_n]\backslash \dd)$ and $ ([k_1]\times  \dots \times [k_n]\backslash \dd)$.
\end{Thm}

\subsection{Proof of the  theorem \ref{T:baw}} \label{S:emo}

\begin{pf}
First we try to  exhibit a simple geometric meaning of a diagram:
\begin{equation}
[k_1]\times  \dots \times [k_n]@>\phi>>\dd(\g)
\end{equation}
Indeed we can think of an object $[k_1]\times  \dots \times [k_n]$ as of a string of polygons. The map $\phi=\phi_1\times  \dots \times \phi_n$ maps the vertices of $[k_i]$ to $i$-th boundary component of $\g$. As a result we have a new structure: a fat graph $\g$ with a collection of labels 
\begin{equation}\label{E:gtz}
b_1(i)\geq \dots \geq b_{k_i}(i)\geq b_1(i)
\end{equation}
 attached to vertices of $i$-th boundary component(and this is done for each boundary component). The above order  in the standard cyclic order, which exists on each boundary.
It is useful to keep in mind the picture \ref{P:lkj}.

  The category   $ ([k_1]\times  \dots \times [k_n]\backslash \dd)$ is a modification of $\tmm$.  The morphisms are  morphisms of $\tmm$ which preserve the labeling. To avoid a possible confusion let us notice  that vertex on the boundary can carry several labels as long the condition \ref {E:gtz} is satisfied.

Observe that we don't need to prove the theorem for all morphisms $\psi$. We can restrict our attention only to generators . Moreover we can disregard  automorphisms $id \times \dots \times \tau \times \dots \times id$ and we may suppose that $\psi=\psi_1\times id \times \dots \times id$.

Let $v:[k_1-1]\rightarrow [k_1]$ be a morphism which inserts a vertex $b_1$ such that $b_{k_1}<b_1<b_2\dots <b_{k_1}$.The morphism $e:[k_1]\rightarrow [k_1-1]$ collapses the edge  connecting $b_1$ and $b_2$.

Observe that $e \circ v=id$. If we can show that the functor  $(v \circ e)^*\overset{def}{=}((v \circ e)\times id \times \dots \times id)^*$ is homotopic to $id$,then $v^*\overset{def}{=}(v\times id \times \dots \times id)^*$ and $e^*\overset{def}{=}(e\times id \times \dots \times id)^*$ would be homotopy inverse functors. This would prove the claim.

Let us examine is the action of the functors $e^*$ and $v^*$ on the categories. Since the map $v$ inserts a vertex $b_1$, the functor $v^*$ remove the label $b_1$ from the set of labels of the boundary component $\partial_1$. The map $e$ collapses the edge between $b_1$ and $b_2$.  the functor $e^*$ adds the label $b_1$ to the vertex which already had label $b_2$.

On the enclosed picture we illustrate the action of the functor $(v \circ e)^*$. It removes the label from the vertex and  adds it to a vertex which has the label $b_2$.

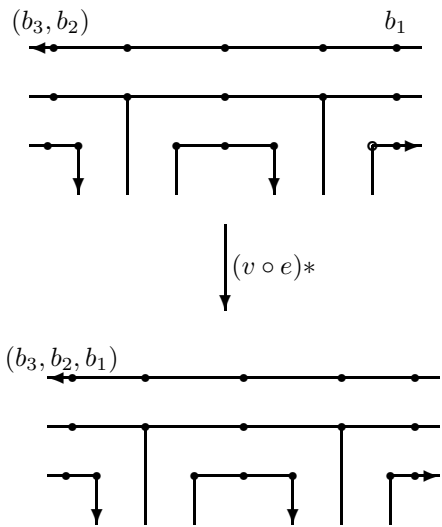
\begin{equation}
\end{equation}

\setlength{\unitlength}{0.00043300in}

\begin{picture}(5337,6291)(-1000,-6223)
\thicklines
\put(2626,-5011){\circle*{100}}
\put(3826,-4411){\circle*{100}}
\put(3826,-5011){\circle*{100}}
\put(3826,-5611){\circle*{100}}
\put(5626,-5611){\circle*{100}}
\put(5026,-5011){\circle*{100}}
\put(4426,-5611){\circle*{100}}
\put(5026,-4411){\circle*{100}}
\put(2026,-5611){\circle*{100}}
\put(2626,-4411){\circle*{100}}
\put(3226,-5611){\circle*{100}}
\put(1651,-5611){\circle*{100}}
\put(1726,-5011){\circle*{100}}
\put(1726,-4411){\circle*{100}}
\put(5926,-5611){\circle*{100}}
\put(5926,-4411){\circle*{100}}
\put(5926,-5011){\circle*{100}}
\put(2401,-961){\circle*{100}}
\put(3601,-361){\circle*{100}}
\put(3601,-961){\circle*{100}}
\put(3601,-1561){\circle*{100}}
\put(5401,-1561){\circle{100}}
\put(4801,-961){\circle*{100}}
\put(4201,-1561){\circle*{100}}
\put(4801,-361){\circle*{100}}
\put(1801,-1561){\circle*{100}}
\put(2401,-361){\circle*{100}}
\put(3001,-1561){\circle*{100}}
\put(1426,-1561){\circle*{100}}
\put(1501,-961){\circle*{100}}
\put(1501,-361){\circle*{100}}
\put(5701,-1561){\circle*{100}}
\put(5701,-361){\circle*{100}}
\put(5701,-961){\circle*{100}}
\put(6226,-4411){\vector(-1, 0){4800}}
\put(1426,-5011){\line( 1, 0){4800}}
\put(2626,-5011){\line( 0,-1){1200}}
\put(5026,-5011){\line( 0,-1){1200}}
\put(1426,-5611){\line( 1, 0){600}}
\put(2026,-5611){\vector( 0,-1){600}}
\put(5626,-6211){\line( 0, 1){600}}
\put(5626,-5611){\vector( 1, 0){600}}
\put(3226,-6211){\line( 0, 1){600}}
\put(3226,-5611){\line( 1, 0){1200}}
\put(4426,-5611){\vector( 0,-1){600}}
\put(6001,-361){\vector(-1, 0){4800}}
\put(1201,-961){\line( 1, 0){4800}}
\put(2401,-961){\line( 0,-1){1200}}
\put(4801,-961){\line( 0,-1){1200}}
\put(1201,-1561){\line( 1, 0){600}}
\put(1801,-1561){\vector( 0,-1){600}}
\put(5401,-2161){\line( 0, 1){600}}
\put(5401,-1561){\vector( 1, 0){600}}
\put(3001,-2161){\line( 0, 1){600}}
\put(3001,-1561){\line( 1, 0){1200}}
\put(4201,-1561){\vector( 0,-1){600}}
\put(3601,-2536){\vector( 0,-1){1050}}
\put(901,-4261){\makebox(0,0)[lb]{$(b_3,b_2,b_1)$}}
\put(5551,-136){\makebox(0,0)[lb]{$b_1$}}
\put(976,-136){\makebox(0,0)[lb]{$(b_3,b_2)$}}
\put(3676,-3136){\makebox(0,0)[lb]{$(v \circ e)*$}}
\end{picture}

It would be nice to have a morphism $id \rightarrow (v \circ e)^*$, which  provided a homotopy. Unfortunately it is not the case.

Instead we shall construct a pair of functors $\Phi$ and $\Psi$ . The first functor satisfy an equation $\Phi^N= (v \circ e)^*$. The power $N$ depends only on the genus $g$ the number of boundary components $n$ and $k_1, \dots, k_n$. The second functor $\Psi$ is equipped with morphisms $\alpha:\Psi\rightarrow id$ and $\beta:\Psi \rightarrow \Phi$. It is clear that this is enough for $id$ and $ (v \circ e)^*$ to be homotopic.

\begin{Def}
Remember that in definition \ref{D:,qe} we identified the vertices of a boundary polygon with ordered triples $(f_1,f_2,v)$. we shall call the valence of  $(f_1,f_2,v)$ the valence of the vertex $v$.
\end{Def}
\begin{Def}
Suppose the vertex $(f_1,f_2,v)$ has a valence two. A vertex\\ $(f_2,f_1,v)$ will be called the adjacent vertex. The reader must  keep in mind that the adjacent vertex may belong to a different boundary component. 
\end{Def}

\begin{Def} of the functor $\Phi$. 
First of all this functor is identical on underlying objects and morphisms of the category $\tmm$.

Fix an object $Y \in Ob([k_1] \times \dots \times [k_n] \backslash \dd )$ . Denote the vertex on the boundary to which $b_1$ is attached by $l$.  If $l$ also carries a label $b_2$, then $\Phi(Y)=Y$. Otherwise choose the nearest staying ahead of $l$ on the same boundary component  vertex $m$ for which one of the below conditions are met:
 
a)A vertex $m$ has a valence two and it or the adjacent vertex carries a label.

b)A vertex $m$ has  a valence grater then two.

    The functor $\Phi$ removes the label $b_1$ from $l$ and puts it on $m$. 
\end{Def}
It should be clear that sufficient high power if this functor gives $(v \circ e)^*$

It is convenient to show how to define the functor $\Psi$ graphically.

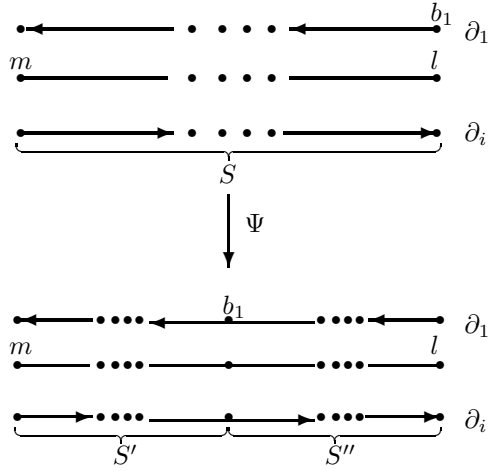
\begin{equation}\label{E:psm}
\end{equation}

\setlength{\unitlength}{0.00043300in}
\begin{picture}(5269,5947)(0,-5489)
\thicklines
\put(1051,-1261){\circle*{100}}
\put(1051,-1936){\circle*{100}}
\put(1051,-661){\circle*{100}}
\put(3526,-661){\circle*{100}}
\put(3826,-661){\circle*{100}}
\put(4126,-661){\circle*{100}}
\put(3151,-661){\circle*{100}}
\put(3526,-1261){\circle*{100}}
\put(3826,-1261){\circle*{100}}
\put(4126,-1261){\circle*{100}}
\put(3151,-1261){\circle*{100}}
\put(3526,-1936){\circle*{100}}
\put(3826,-1936){\circle*{100}}
\put(4126,-1936){\circle*{100}}
\put(3151,-1936){\circle*{100}}
\put(6151,-661){\circle*{100}}
\put(6151,-1261){\circle*{100}}
\put(6151,-1936){\circle*{100}}
\put(1013,-4786){\circle*{100}}
\put(1013,-5413){\circle*{100}}
\put(1013,-4229){\circle*{100}}
\put(2215,-4229){\circle*{100}}
\put(2362,-4229){\circle*{100}}
\put(2507,-4229){\circle*{100}}
\put(2033,-4229){\circle*{100}}
\put(2215,-4786){\circle*{100}}
\put(2362,-4786){\circle*{100}}
\put(2507,-4786){\circle*{100}}
\put(2033,-4786){\circle*{100}}
\put(2215,-5413){\circle*{100}}
\put(2362,-5413){\circle*{100}}
\put(2507,-5413){\circle*{100}}
\put(2033,-5413){\circle*{100}}
\put(4915,-4229){\circle*{100}}
\put(5062,-4229){\circle*{100}}
\put(5207,-4229){\circle*{100}}
\put(4733,-4229){\circle*{100}}
\put(4915,-4786){\circle*{100}}
\put(5062,-4786){\circle*{100}}
\put(5207,-4786){\circle*{100}}
\put(4733,-4786){\circle*{100}}
\put(4915,-5413){\circle*{100}}
\put(5062,-5413){\circle*{100}}
\put(5207,-5413){\circle*{100}}
\put(4733,-5413){\circle*{100}}
\put(6192,-4229){\circle*{100}}
\put(6192,-4786){\circle*{100}}
\put(6192,-5413){\circle*{100}}
\put(3601,-5413){\circle*{100}}
\put(3601,-4786){\circle*{100}}
\put(3601,-4229){\circle*{100}}
%\put(2251,314){\vector(-2,-1){1020}}
\put(3601,-2686){\vector( 0,-1){900}}
\put(1051,-1261){\line( 1, 0){1800}}
\put(4351,-1261){\line( 1, 0){1800}}
\put(2926,-661){\vector(-1, 0){1800}}
\put(6151,-661){\vector(-1, 0){1800}}
\put(1051,-1936){\vector( 1, 0){1875}}
\put(4276,-1936){\vector( 1, 0){1875}}
\put(1013,-4786){\line( 1, 0){875}}
\put(1924,-4229){\vector(-1, 0){875}}
\put(1013,-5413){\vector( 1, 0){911}}
\put(5317,-4786){\line( 1, 0){875}}
\put(6192,-4229){\vector(-1, 0){875}}
\put(6492,-790){\makebox(0,0)[lb]{$\partial _1$}}
\put(6492,-2040){\makebox(0,0)[lb]{$\partial _i$}}

\put(6492,-4390){\makebox(0,0)[lb]{$\partial _1$}}
\put(6492,-5560){\makebox(0,0)[lb]{$\partial _i$}}
\put(5280,-5413){\vector( 1, 0){912}}
\put(2626,-5461){\vector( 1, 0){2025}}
\put(2626,-4786){\line( 1, 0){2025}}
\put(4651,-4261){\vector(-1, 0){2025}}
%\put(2401,314){\makebox(1000,1000)[lb]{$\Phi(b_1)$}}
\put(3801,-3150){\makebox(0,0)[lb]{$\Psi$}}
\put(6076,-561){\makebox(0,0)[lb]{$b_1$}}
\put(6076,-1150){\makebox(0,0)[lb]{$l$}}
\put(913,-1150){\makebox(0,0)[lb]{$m$}}
\put(6076,-4650){\makebox(0,0)[lb]{$l$}}
\put(913,-4650){\makebox(0,0)[lb]{$m$}}
\put(3480,-2550){\makebox(0,0)[lb]{$S$}}
\put(2180,-6000){\makebox(0,0)[lb]{$S'$}}
\put(4780,-6000){\makebox(0,0)[lb]{$S''$}}

\put(3526,-4161){\makebox(0,0)[lb]{$b_1$}}
\thinlines
\put(2180,-5550){\oval(2400.0,200.0)[bl]}
\put(2380,-5550){\oval(2400.0,200.0)[br]}
\put(2380,-5750){\oval(200.0,200.0)[tl]}
\put(2180,-5750){\oval(200.0,200.0)[tr]}

\put(4800,-5550){\oval(2400.0,200.0)[bl]}
\put(5000,-5550){\oval(2400.0,200.0)[br]}
\put(5000,-5750){\oval(200.0,200.0)[tl]}
\put(4800,-5750){\oval(200.0,200.0)[tr]}

\put(3500,-2070){\oval(5000.0,200.0)[bl]}
\put(3700,-2070){\oval(5000.0,200.0)[br]}
\put(3700,-2270){\oval(200.0,200.0)[tl]}
\put(3500,-2270){\oval(200.0,200.0)[tr]}
\end{picture}
\begin{equation}\notag
\end{equation}
On the picture \ref{E:psm} the reader  may see a segment $S$ of a graph $\g$ which is being transformed by the functor $\Psi$. The end points of this segment are $l$ and $m$ , which where used in the definition of functor $\Phi$.The segment might be a union of a several edges. The functor $\Psi$ excises the segment $S$, makes two copies of it which we denote by $S'$ and $S''$, glues them and  pastes the whole thing  back. The point $b_1$ is placed  in the way the reader can see on the picture.

It is easy to see that $\Psi$ is a functor. The desired morphisms are
$S'':\Psi \rightarrow Id$ (it collapses the segment $S''$), $S':\Psi \rightarrow \Phi$(it collapses the segment $S'$).     
\end{pf}
\begin{Prop}
The nerve of the category $([1]\times \dots \times[1]\backslash \dd)$ is homotopic to $K(\g_{g,n},1)$. The group $\g_{g,n}$ is the mapping class group of orientation preserving diffeomorphisms of a surface $X$ with $n$ boundary components, acting trivially on the boundary. 
\end{Prop}
\begin{pf}
We start with some simple reductions. Let $P:([1]\times \dots \times[1]\backslash \dd)\rightarrow \tmm$ be the forgetting functor. The category $\mm$ sits inside of $\tmm$ as a full subcategory (see proposition \ref {P:wbi}). Consider the inverse image $i:P^{-1}(\mm)\hookrightarrow ([1]\times \dots \times[1]\backslash \dd)$. Define a homotopy inverse functor $er: ([1]\times \dots \times[1]\backslash \dd) \rightarrow P^{-1}(\mm)$. It erase all vertices of valence two with and move their labels to the next ahead vertex of valence grater or equal to three.
The arguments that $er$ and $i$ are homotopy inverse follow closely the proof of the previous theorem.  

We want to introduce an alternative definition of the category $P^{-1}(\mm)\overset{def}{=}\mn$ which will prove to be more convenient for our purposes.

 In this new definition the objects of $\mn$ will be certain fat graphs whose  vertices may have  valence one . The restriction which distinguish a graph as being an object of $\mn$ is each orbit of the  transformation $D:fl(\g) \rightarrow fl(\g)$ contains exactly one flag adjacent to a vertex of valence one, one-valent vertices are connected to vertices of valence four or grater and the graph has no two-valent vertices. The morphisms are defined as before except we do not allow to shrink an edge whose one end is one-valent vertex.

The equivalence of these two definitions follows from existence of a  functor which maps  graphs $\g \in \mm$ with one label  on each  boundary component to a graphs  with one-valent vertices. If a label $b$ is attached to a vertex $l$ the functor  introduces a new vertex $\tilde b$ which is  connected with $l$ by an edge $e$.

 The flag which corresponds to the $l$-end of $e$ we denote by $f_e$.  By definition \ref{D:,qe}  we identify the label $b$ with an ordered triple $(f_1,f_2,l)$ where $f_2$ is  a flag staying ahead of $f_1$ in the cyclic order of the vertex $l$. The cyclic order on the new graph will be $f_1<f_e<f_2<f_3 \dots$.
This identification enables us to define a boundary  component $\partial(v)$ corresponding to a one-valent vertex $v$.

Fix an oriented surface $X$ with a boundary. Assume that on  each boundary component $\partial_i(X)$ a point $a_i$ is chosen.
Define a functor $f:\mn \rightarrow Set$. The set $f(\g)$ consists of isotopy classes of coherent embeddings of $CW(\g)$ such that the complement to one-valent vertices sits in the complement to the boundary. The image of each one-valent vertex $v$ is $a_i$, the index $i$ is the one of $\partial_i(v)$.

Arguing as in the proof of the proposition \ref{P:aie} one can show that the action of the  group $\g_{g,n}$ is free and transitive on $f(\g)$. It is also evident that the functor $f$ transforms morphisms into isomorphisms, therefore the category $(\mn,f)$  is a covering of $\mn$.
Now everything follows from the 
\begin{Thm}\cite{H}
The category $(\mn,f)$ is contractible.
\end{Thm}
\end{pf}
\begin{Cor}
The homology of the category $\mn$ are zero above dimension $6g-6+3n$
\end{Cor}
\begin{pf}
The classifying space $K(\g_{g,n},1)$ is homotopy equivalent to the moduli of curves with $n$ marked points and unit tangent vectors to these points. This is a smooth manifold of the  dimension $6g-6+3n$.
\end{pf}

\subsection{A spectral sequence}

\begin{Prop}\label {P:wmi}
For any functor $G:\Lambda \times \dots \times \Lambda \rightarrow Vect$ there is a spectral sequence $E^k_{.,.}$ which converges to $H_.(\tmm,(\dd)^*G)$, whose second term is equal to $H_.(\Lambda \times \dots \times \Lambda,G)\otimes H_.(\g_g^n)$. The spectral sequence degenerates in $6g-6+3n$ term.
\end{Prop}
\begin{pf}
We take advantage of the spectral sequence from section \ref{S:nnb}. According to  the theorem \ref{T:baw} the homotopy type of the category $ ([k_1]\times  \dots \times [k_n]\backslash \dd)$ is independent on $[k_1]\times  \dots \times [k_n]$ and equal to $K(\g_{g,n},1)$ Due to remark \ref{R:ngc} we can argue that the second term of our spectral sequence is precisely equal to the one described in the statement of the proposition.
\end{pf}

\section{Example of computation of graph homology.}

 To formulate a statement we need to introduce a category $\D$. Its objects are sets $[n]=\{1 <2< \dots <n\}$ with the specified linear order. The morphisms are nondecreasing maps. The category $\D$ contains a subcategory $\D_2$ . The sets $[n]$ ($n>1$) make the objects. A morphism $\alpha:[n]\rightarrow [k]$ is a morphism in  $\D$ which map the minimal element of $[n]$ to the minimal element of $[k]$, similarly with maximal elements. It is a standard exercise to check that $\D$ and $\D_2$ are antiisomorphic. On the objects the antiisomorphism  $i:\D \rightarrow \D_2$ maps $[n]$ to $[n+1]$. On the morphisms it maps inclusions into surjections and  surjections into inclusions . As a result any covariant functor on $\D_2$ can be thought of as a contravariant functor on $\D$, such is usually called a simplicial object. In particular if $G: \D_2 \rightarrow Vect$ a covariant functor to vector spaces, then $H_.(\D_2,G)=H_.(G)$. The right group is a ordinary homology of a simplicial vector space, which can be computed via standard complex.
\begin{Rem}\label{R:wjc} The categories  $([1]\backslash Id_{\Lambda})$ and  $\D_2$ are isomorphic. As we know a  diagram $[1] \rightarrow [n] \in Ob(([1]\backslash Id_{\Lambda}))$ can be thought  as a cyclic graph with a choice of a label on one of it vertices. This label specifies a minimal and a maximal  element in the order which exists on vertices of $[n]$ . Therefore we identified objects of $([1]\backslash Id_{\Lambda})$ with ones of $\D_2$. It is clear that this identification is functorial.
\end{Rem}
There is a projection $p:([1]\backslash Id_{\Lambda}) \rightarrow \Lambda$. 
\begin{Prop}\label {P:awj}
For any covariant functor $G:\Ltr \rightarrow Vect$ there is an isomorphism $H_.(\mtr,(\dtr)^*G)=H_.(\D_2 \times \D_2 \times \D_2,(p_1 \circ i_1 \times p_2 \circ i_2 \times p_3 \circ i_3 )^*G)$
\end{Prop}
\begin{pf}     
We want to address the following question first. Are  there  functors 
$\alpha_i$($i=1,2,3$) which would make the diagram below to commute?
\begin{equation}
\begin{CD}
\mtr  @>\text{$\partial_i$}>> \Lambda\\
@VV\text{$\alpha_i$}V   @AA\text{$p_i$}A \\
([1]\backslash \partial_i)  @>\text{$\tilde \partial_i$}>> ([1]\backslash Id_{\Lambda})
\end{CD} 
\end{equation}
The reason why we anticipate that such functor might exist is that the nerve of $\mtr$ is contractible, therefore there no obstruction to existence of appropriate map of nerves. We want to construct a functor which would realize such a map.
 Actually we shall construct only $\alpha_1$. Two other functors are defined along the same lines. As it was explained in the beginning of section \ref{S:emo} the objects of $([1]\backslash \partial_i)$ are fat graphs from $\mtr$ with a choice of one label on the first boundary component. We expect that the functor $\alpha_1$ only places a label in the appropriate spot on the boundary, but does nothing with the underlying graph. Below we directly exhibit all graphs from $\mtr$ which have no two-valent vertices, and spots where the functor places the label:

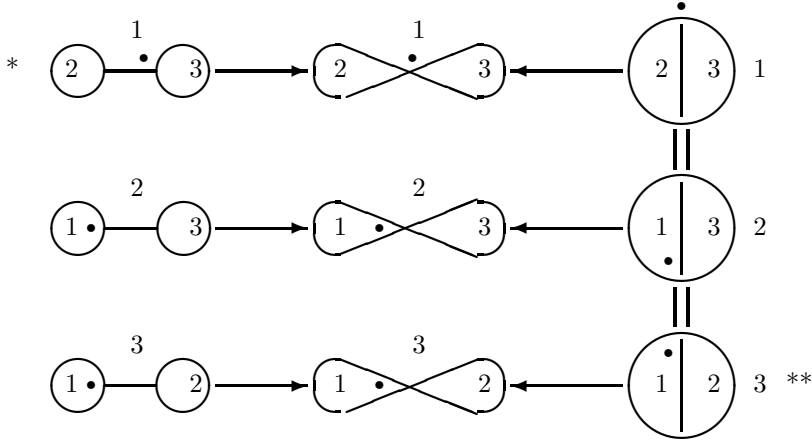
\setlength{\unitlength}{0.00045300in}

\begin{picture}(8925,5045)(376,-4573)
\thicklines
\put(4201,-2161){\oval(600,600)[tl]}
\put(4201,-2161){\oval(600,600)[bl]}
\put(5776,-2161){\oval(600,600)[br]}
\put(5776,-2161){\oval(600,600)[tr]}
\put(4201,-3961){\oval(600,600)[tl]}
\put(4201,-3961){\oval(600,600)[bl]}
\put(5776,-3961){\oval(600,600)[br]}
\put(5776,-3961){\oval(600,600)[tr]}
\put(4201,-361){\oval(600,600)[tl]}
\put(4201,-361){\oval(600,600)[bl]}
\put(5776,-361){\oval(600,600)[br]}
\put(5776,-361){\oval(600,600)[tr]}
\put(1201,-361){\circle{600}}
\put(1201,-2161){\circle{600}}
\put(1201,-3961){\circle{600}}
\put(2401,-361){\circle{600}}
\put(2401,-3961){\circle{600}}
\put(2410,-2161){\circle{600}}
\put(8101,-361){\circle{1200}}
\put(8101,-2161){\circle{1200}}
\put(8101,-3961){\circle{1200}}
\put(1951,-211){\circle*{100}}
\put(1351,-2161){\circle*{100}}
\put(1351,-3961){\circle*{100}}
\put(5026,-211){\circle*{100}}
\put(4651,-2161){\circle*{100}}
\put(4651,-3961){\circle*{100}}
\put(7951,-3586){\circle*{100}}
\put(7951,-2536){\circle*{100}}
\put(8101,389){\circle*{100}}
\put(1501,-3961){\line( 1, 0){600}}
\put(1501,-2161){\line( 1, 0){600}}
\put(1501,-361){\line( 1, 0){600}}
\put(4201,-2461){\line( 5, 2){1564.655}}
\put(4201,-1861){\line( 5,-2){1564.655}}
\put(8101,170){\line( 0,-1){1050}}
\put(8101,-1640){\line( 0,-1){1050}}
\put(8101,-3440){\line( 0,-1){1050}}
\put(2776,-3961){\vector( 1, 0){1050}}
\put(2776,-2161){\vector( 1, 0){1050}}
\put(2776,-361){\vector( 1, 0){1050}}
\put(7426,-3961){\vector(-1, 0){1275}}
\put(7426,-2161){\vector(-1, 0){1275}}
\put(7426,-361){\vector(-1, 0){1275}}
\put(4276,-4261){\line( 5, 2){1500}}
\put(4201,-3661){\line( 5,-2){1564.655}}
\put(4276,-661){\line( 5, 2){1500}}
\put(4201,-61){\line( 5,-2){1564.655}}
\put(8026,-1036){\line( 0,-1){450}}
\put(8176,-1036){\line( 0,-1){450}}
\put(8176,-2836){\line( 0,-1){450}}
\put(8026,-2836){\line( 0,-1){450}}
\put(1051,-4036){\makebox(0,0)[lb]{1}} 
\put(2476,-4036){\makebox(0,0)[lb]{2}} 
\put(5776,-4036){\makebox(0,0)[lb]{2}} 
\put(7801,-4036){\makebox(0,0)[lb]{1}} 
\put(1801, 14){\makebox(0,0)[lb]{1}} 
\put(5026, 14){\makebox(0,0)[lb]{1}} 
\put(1801,-1786){\makebox(0,0)[lb]{2}} 
\put(5026,-1786){\makebox(0,0)[lb]{2}} 
\put(1801,-3586){\makebox(0,0)[lb]{3}} 
\put(5026,-3586){\makebox(0,0)[lb]{3}} 
\put(1051,-2236){\makebox(0,0)[lb]{1}} 
\put(1051,-436){\makebox(0,0)[lb]{2}} 
\put(2476,-436){\makebox(0,0)[lb]{3}} 
\put(2476,-2236){\makebox(0,0)[lb]{3}} 
\put(4126,-4036){\makebox(0,0)[lb]{1}} 
\put(4126,-2236){\makebox(0,0)[lb]{1}} 
\put(4126,-436){\makebox(0,0)[lb]{2}} 
\put(5776,-2236){\makebox(0,0)[lb]{3}} 
\put(5776,-436){\makebox(0,0)[lb]{3}} 
\put(7801,-2236){\makebox(0,0)[lb]{1}} 
\put(7801,-436){\makebox(0,0)[lb]{2}} 
\put(8401,-4036){\makebox(0,0)[lb]{2}} 
\put(8401,-2236){\makebox(0,0)[lb]{3}} 
\put(8401,-436){\makebox(0,0)[lb]{3}} 
\put(8926,-4036){\makebox(0,0)[lb]{3}} 
\put(8926,-2236){\makebox(0,0)[lb]{2}} 
\put(8926,-436){\makebox(0,0)[lb]{1}} 
\put(376,-436){\makebox(0,0)[lb]{*}} 
\put(9301,-4036){\makebox(0,0)[lb]{**}} 
\end{picture}

Two  comments about the picture are in order. The functor $\alpha_1$ is defined not complete canonically. On the picture we marked objects (*) and (**) where we have a choice of placing the label. The equality signs mean that the objects are isomorphic. The arrows denote collapsing of appropriate edges. It is clear from the picture that the rule of placing labels is manifestly functorial.On the objects with two-valent vertices we extend the above rule by functoriality. 
It follows that $H_.(\mtr,(\dtr)^*G)=H_.(\mtr,(\ptr \circ \tilde \dtr \circ \alphatr)^*G)$.

\begin{Lem}\label{L:heb}
For any object $Y \in ([1] \times [1]\times [1] \backslash \dtr)$ the category $(Y\backslash \alphatr$ is contractible
\end{Lem}
\begin{Lem}\label{L:pna}
For any object $X \in ([1] \times [1] \times [1] \backslash \ptr)$ the category $(Y\backslash \tilde\dtr)$ is contractible
\end{Lem}
We don't give proves since they don't use any new ideas.

 We split computation of $H_.(\mtr,(\ptr \circ \tilde \dtr\circ \alphatr)^*G)$ in two steps.
Utilizing a spectral sequence described in  \ref{S:nnb} applied to $(\ptr)^*G$ and a functor $\tilde \dtr$, we see that due to lemma \ref{L:pna} the spectral sequence  degenerates in the first term. We deduce the identity 
\begin{align}
H_.([1]&\times[1]\times[1] \backslash Id_{\Ltr},(\ptr)^*G)=\notag\\
&H_.([1]\times[1]\times[1] \backslash \tilde \dtr ,(\ptr\circ \tilde \dtr)^*G)\notag
\end{align}

We repeat the argument with spectral sequence one more time but use as input  $(\ptr\circ \tilde \dtr)^*G$ and the functor $\alphatr$. This time we use lemma \ref{L:heb}. We see again the the groups $H_.(([1]\times[1]\times[1] \backslash \tilde \dtr ,(\ptr\circ \tilde \dtr)^*G)$ and $H_.(\mtr ,(\ptr\circ \tilde \dtr \circ \alphatr)^*G)$ are equal.

So we have 
\begin{align}
H_.(\mtr& ,(\ptr\circ \tilde \dtr \circ \alphatr)^*G)=\notag\\
&H_.(([1]\times[1]\times[1] \backslash Id_{\Ltr},(\ptr)^*G)\notag
\end{align}
According to remark \ref{R:wjc} the later groups  are just homology groups of simplicial vector space associated with  $G$. 
\end{pf} 

\section{Application to algebras}
Fix an associative algebra $A$ with a unit and define a functor $A^{\#}$ from the category $\Lambda$ to vector spaces. The objects of $\Lambda$ are $n$-gons. The functor assigns to an $n$-gone a vector space equal to $A^{\otimes n}$. In the space $A^{\otimes n}$ each tensor factor correspond to a vertex of the $n$-gone. If we think of an $n$ -gone as of a cyclic graph the functor $A^{\#}$ transforms collapsing of an edge into multiplication of the tensor factors corresponding to the ends of the edge. The order in which multiplication is carried out is  determined by the orientation of the cyclic graph. The group $S_n$ acts on the space   $A^{\otimes n}$ by permutations, the action of rotation of the $n$-gone factors through the action of $S_n$. The functor  $A^{\#}$ transforms insertions of a vertex into an edge into tensoring by the unit , where the unite is placed exactly on the spot of the tensor factor corresponding to the new vertex.

We also give a definition of such functor for non unital algebras. Given an algebra $I$ with no unit we define a new algebra $I^+$ with adjoint unit. The algebra $I^+$ admits an augmentation  $p:I^+ \rightarrow k$, where $k$ is the ground field. Therefore we have a map $p:I^{+\#}  \rightarrow k^{\#}$. Denote $I^{\#}$ the kernel of this map.

\begin{Def} Fix $n$ algebras with a unit. Denote 1111
\begin{equation}
H_.^{g}(A_1,\dots,A_n)= \underset{\overset{\longrightarrow}{\tmm}}{\lim}^.  \bigotimes_{i=1}^n A_i^{\#} \circ\partial_i
\end{equation}
the graph homology of algebras.
\end{Def}
\begin{Rem}
For the case on non-unital algebras we replace the functor $A^{\#}$ by $I^{\#}$.
\end{Rem}

\begin{Rem}
One can replace the category $\tmm$ by category $\tm$. However they compute the same homology groups as was shown in proposition \ref {P:ghe}.
\end{Rem}

\begin{Rem}
According to Connes \cite{C} for one algebra the groups $ \underset{\overset{\longrightarrow}{\Lambda}}{\lim}^.A^{\#}$ are equal to cyclic homology $HC_.(A)$. It explains why  graph homology is a generalization of cyclic homology.
\end{Rem}

\begin{Rem}
Due to identification ${\Bbb T'M}_{0,2}=\Lambda$ we have $H_.^{0}(A_1,A_2)=HC_.(A_1 \otimes A_2^{op})$.
\end{Rem}

The proposition \ref{P:wmi} in context of functor $G=\bigotimes_{i=1}^nA_i^{\#} \circ\partial_i$ becomes 
\begin{Prop}
There is a spectral sequence converging to $H_.^{g}(A_1,\dots,A_n)$
whose second term is equal to $(\bigotimes_{s=1}^nHC_.(A_s))_i\otimes H_j(\g_{g,n},k)$
\end{Prop}
\begin{Cor}
Graph homology are Morita invariant
\end{Cor}
\begin{pf}
We establish it in the following case. Suppose maps  $A_i \rightarrow A'_i$ are inclusions of Morita equivalent algebras. We show that these homomorphisms induce isomorphism on graph homology.

The homomorphisms define a morphism between spectral sequences of graph homology of algebras $A_1, \dots, A_n$ and  $A'_1, \dots, A'_n$ . It is known (see \cite {L}) that the maps $HC_.(A_i) \rightarrow HC_.(A'_i)$ on cyclic homology are  isomorphisms. Therefore the second terms of the spectral sequences of graph homology are isomorphic. It implies that  limit terms must be isomorphic too.   
\end{pf}
\begin{Cor} Now we can make a general computation of zero graph homology. Indeed from the spectral sequence it follow that $H^g_0(A_1, \dots,A_n)=H_0(A_1)\otimes \dots \otimes H_0(A_n)$ - the tensor product of zero  Hochschild homology.
\end{Cor}

The proposition \ref {P:awj} specializes in our case into 
\begin{Prop}
There is an isomorphism $H^0_.(A_1,A_2,A_3)=H_.(A_1)\otimes H_.(A_2)\otimes H_.(A_3)$. The groups $H_.(A)$ are Hochschild homology groups of an algebra $A$.
\end{Prop}

Let us examen the relative case of graph homology. Suppose we have an ideal $I \subset A$. Denote by $p$ the natural projection $A \rightarrow A/I$. It induces a map $p^{\#}:A^{\#} \rightarrow (A/I)^{\#}$. Denote the  kernel of $p^{\#}$  by $(A,I)^{\#}$. There is a natural map $I^{\#} \rightarrow (A,I)^{\#}$.

\begin{Def}
Define a relative graph homology by the formula 
\begin{align}
H_.^{g}&(A_1,\dots,(A_i,I_i),\dots,A_n)=\notag\\
& \underset{\overset{\longrightarrow}{\tmm}}{\lim}^. A_1^{\#} \circ\partial_1\otimes,\dots, \otimes (A_i,I_i)^{\#} \circ\partial_i \otimes,\dots, \otimes  A_n^{\#} \circ\partial_n \notag
\end{align}
\end{Def}
By functoriality we have a map
\begin{equation}\label{E:vae}
H_.^{g}(A_1,\dots,I_i,\dots,A_n) \rightarrow H_.^{g}(A_1,\dots,(A_i,I_i),\dots,A_n)
\end{equation}
We say that an ideal is excisive in graph homology if the map \ref {E:vae} is an isomorphisms.
\begin{Prop}
If an ideal $I$ is excisive for cyclic homology then it is excisive in graph homology.
\end{Prop}
\begin{pf}
We take advantage of our spectral sequence. According to the  assumptions the second terms of the spectral sequences for $ H_.^{g}(A_1,\dots,I_i,\dots,A_n)$ and $H_.^{g}(A_1,\dots,(A_i,I_i),\dots,A_n)$ are isomorphic. Therefore the limiting terms must coincide. 
\end{pf}
\begin{Cor}
As usual for excisive ideal we have a long exact sequence of graph homology :
\begin{align}
&\rightarrow  H_s^{g}(\dots,I_i,\dots)\rightarrow H_s^{g}(\dots,A_i,\dots)\rightarrow \\
&\rightarrow H_s^{g}(\dots,(A_i/I_i),\dots)\rightarrow  H_{s-1}^{g}(\dots,I_i,\dots)\rightarrow
\end{align}
\end{Cor}

\section{A definition of graph K-theory}
 Denote by  $\m{}{A}$ the  algebra of $r
\times r$ matrices
with entries in algebra $A$. we shall suppress some symbols in the
notations for graph homology of the ground  field,i.e:
\begin{equation}
H_.^g(n,\k) \overset{def}{=}H_.^g(\k,\dots,\k)
\end{equation}

We are going to introduce a class of   morphisms :
\begin{equation}
H_.^g(n,\k) \rightarrow H_.^g(\m{1}{A_1},\dots,\m{n}{A_{n}})
\end{equation}

we shall do it on the  chain level. Let us invoke a functor

\begin{equation}
L
\begin{Sb}
{\bold Mat}
\begin{Sb}
r_1
\end{Sb}(A_1)
\dots
{\bold Mat}
\begin{Sb}
r_n
\end{Sb}(A_n)
\end{Sb}
\overset{def}{=} \bigotimes_{i=1}^n \m{i}{A_1}^{\#} \circ \partial_i: \tm \to Vect
\end{equation}
used in the definition of graph homology.  We use  a notation $L(\Gamma)$ for a degenerate case
when all algebras $A_1= \dots A_n= \k$. In this case $L(\Gamma)=\k$ for
all $\Gamma$

Suppose also  $\psi(\Gamma)$ is a morphism of functors, i.e.
\begin{equation} \label{E:rfy}
\psi(\Gamma):L(\g) \rightarrow \gl{ \g}
\end{equation}
 is a homomorphism of
$Aut(\Gamma)$ modules. 

and all diagrams\\
\begin{equation} \label{E:giy}
\begin{CD}
L(\Gamma)@>L(e)>>L(\Gamma/e)\\
@V\psi(\Gamma)VV  @V\psi(\Gamma/e)VV\\
\gl{\Gamma}@>L(e)>>\gl{\Gamma/e}\\
\end{CD}
\end{equation}
\\
commute.
Then the map $\psi$ induces a morphism
\begin{equation}
H_.^g(n,\k ,\dots,\k)\overset{\psi}{\rightarrow} H_.^g(\m{1}{A_1},\dots,\m{n}{A_{n}})
\end{equation}
We use Morita invariance of graph homology.

The spaces $L(\Gamma)$ are canonically isomorphic to $\k$. The image
of the unit under the  map $\psi$ defines an element in
$\gl{\Gamma}$, denoted by the same letter $\psi(\Gamma)$. 

\begin{Rem} \label{R:vqw}
The set of all these elements are connected by some
relations.For example $\psi(\Gamma/e)=L(e)\psi(\Gamma)$. In particular we can
recover $\psi$ if we know $\psi(\Gamma)$ for graphs $\Gamma$ having  only
vertices of valence two or three. 
\end{Rem}

\subsection{A method to produce $\psi(\Gamma)$}

Fix a vertex $m$ of a fat graph $\g$. Denote $n$ a germ of a boundary curve $s_i$ through $m$ and $N_m$ the set of all germes of various $s_i$ through $m$. Warning: the set $N_m$ might contain more then one germ of the same boundary curve.
\begin{Rem} \label{R:frt}
One can easily see that $\#\{N_m\}=\#\text{flags adjacent to $m$}$.
\end{Rem}

Given a graph $\g$ and its vertex $m$ fix an element
 \begin{equation} \label{E:vjw}
\psi_m(\g) \in \bigotimes
\begin{Sb}
n \  germ \ of s_i \\
n \in N_m
\end{Sb}
 \m{i}{A_i}
\end{equation}

Due to remark \ref{R:frt} 
\begin{equation}
\psi( \Gamma) \overset {def} {=} 
\bigotimes
\begin{Sb}
m \in  vertices \\
of  \g 
\end{Sb} 
\psi_m( \g) \in \gl{ \g}
\end{equation}

Of course in general  collection $\psi(\g)$ will violate conditions \ref {E:rfy}, \ref{E:giy}. According to  remark \ref {R:vqw} it is suffice to fix collections of tensors only for graphs $\g$ with vertices valence two and three.
\begin{Def}\label {D:nsk}
Denote $M_{g,.n}^{2,3}$ the following collection of fat graph. A graph $\g \in M_{g,.n}^{2,3}$ has vertices valence two and three. $\g$ has no two  vertices valence three connected by an edge, similarly $\g$ has no two vertices valence two connected by an edge. In this definition the vertices might coincide.
\end{Def}

 We shall define  a set of equations , whose solutions are $\psi_m(\g)$, $\g \in M_{g,.n}^{2,3}$.

{\bf Equations for quadratic tensors.}
The first set of equations is for $\psi_m(\g)$, where $m$ is a vertex of valence two. We know that according to \ref{E:vjw}
\begin{equation}
\psi_m(\g)=\sum a_l \otimes b_l  \in \m{i}{A_i} \otimes \m{j}{A_j}
\end{equation}
We write the diagram of orientations of germs $p$ $\tilde p$ of $s_i$ and $q$, $\tilde q$ of $s_j$ at a vertex $m$. Suppose it looks like  the upper part of the  following picture:

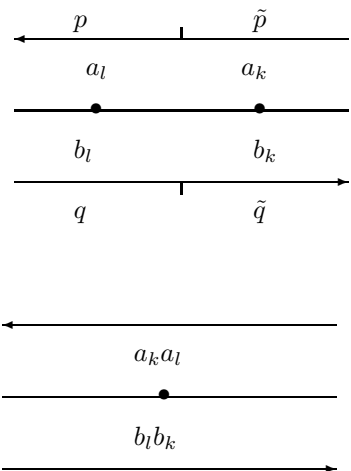
\setlength{\unitlength}{3947sp}%
\begin{picture}(2799,2877)(-100,-2698)%%(2199,2877)(814,-2698)
\thinlines
\put(1373,-460){$\bullet $}
\put(2401,-460){$\bullet $}
\put(1801,-2260){$\bullet $}
\put(901,-436){\line( 1, 0){2100}}
\put(901,-886){\vector( 1, 0){2100}}
\put(2926,-1786){\vector(-1, 0){2100}}
\put(826,-2236){\line( 1, 0){2100}}
\put(826,-2686){\vector( 1, 0){2100}}
\put(1951, 89){\line( 0,-1){ 75}}
\put(1951,-886){\line( 0,-1){ 75}}
\put(3001, 14){\vector(-1, 0){2100}}
\put(2326,-211){\makebox(0,0)[lb]{$a_k$}}
\put(2401,-736){\makebox(0,0)[lb]{$b_k$}}
\put(1351,-211){\makebox(0,0)[lb]{$a_l$}}
\put(1276,-736){\makebox(0,0)[lb]{$b_l$}}
\put(1651,-2536){\makebox(0,0)[lb]{$b_lb_k$}}
\put(1651,-2011){\makebox(0,0)[lb]{$a_ka_l$}}
\put(1276, 89){\makebox(0,0)[lb]{$p$}}
\put(2401, 89){\makebox(0,0)[lb]{$\tilde p$}}
\put(1276,-1111){\makebox(0,0)[lb]{$q$}}
\put(2401,-1111){\makebox(0,0)[lb]{$\tilde q$}}
\end{picture}

 The first set of equations one can read off from the lower part of the picture :
\begin{equation} \label {E:erg}
\sum a_k a_l \otimes b_l b_k=\sum a_l \otimes b_l 
\end{equation}
 Introduce a notation :  suppose $A$ is an algebra, then $A^{op}$ is the algebra with opposite multiplication. By definition it has the same linear space as $A$. Using this notation we can rewrite  \ref{E:erg} shortly
\begin{equation} \label{E:iwf}
\psi_m(\g)^2=\psi_m(\g) \in \m{i}{A_i}^{op} \otimes \m{j}{A_j}
\end{equation}

Notice that if the arrows in the diagram  were  reversed ,  the other tensor multiple ($\m{j}{A_j}$) would acquire  the   $op$ symbol.

Now we shall discuss the meaning of this relation. Suppose we have  a graph $\g'$, such that $\g'/e=\g$ and the ends of the edge  $e$ are vertices  $m$ and $m'$. . We can define $\psi_{ m}(\g')=\psi_{ m'}(\g')\overset{def}{=}\psi_{n}(\g)$, where vertex $n \in \g$ is obtained from $e$ by collapsing it to a point . Then it easy to check that  $L(e)\psi(\g')=\psi(\g)$

{\bf Equations for two-tensors and three-tensors}

Fix vertex $m(2)$ of valence two and $m(3)$ of valence three in a fat graph $\g$. Suppose also they are ends of an edge $e$. Suppose we have the following diagram of orientations of germs of boundary curves at these vertices(the upper part of the picture):

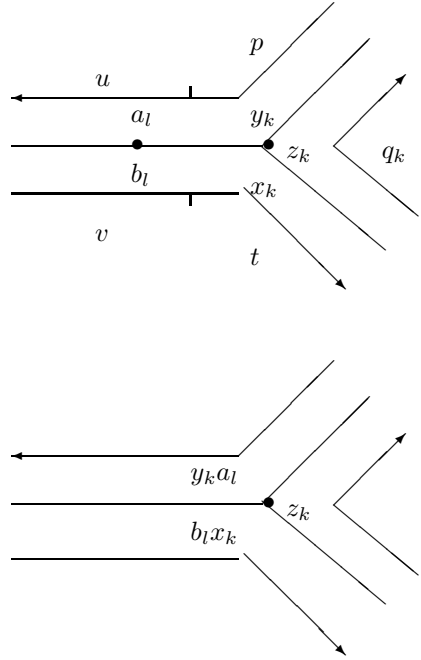
\setlength{\unitlength}{3947sp}%
\begin{picture}(2574,4121)(-200,-4095)%%(2574,4121)
\thinlines
\put(1276,-916){$\bullet $}
\put(2101,-916){$\bullet $}
\put(2101,-3166){$\bullet $}
\put(1988,-1148){\vector( 1,-1){637.500}}
\put(1951,-1186){\line(-1, 0){1425}}
\put(1651,-1186){\line( 0,-1){ 75}}
\put(1651,-511){\line( 0,-1){ 75}}
\put(1951,-586){\vector(-1, 0){1425}}
\put(1951,-586){\line( 1, 1){600}}
\put(2551,-886){\vector( 1, 1){450}}
\put(2551,-886){\line( 6,-5){531.148}}
\put(526,-886){\line( 1, 0){1575}}
\put(2101,-886){\line( 1, 1){675}}
\put(2101,-893){\line( 6,-5){774.590}}
\put(526,-3136){\line( 1, 0){1575}}
\put(2101,-3118){\line( 6,-5){774.590}}
\put(2101,-3136){\line( 1, 1){675}}
\put(2551,-3142){\vector( 1, 1){450}}
\put(2551,-3142){\line( 6,-5){531.148}}
\put(1951,-2836){\vector(-1, 0){1425}}
\put(1951,-2836){\line( 1, 1){600}}
\put(1988,-3445){\vector( 1,-1){637.500}}
\put(1951,-3483){\line(-1, 0){1425}}
\put(2026,-286){\makebox(0,0)[lb]{$p$}}
\put(1051,-1486){\makebox(0,0)[lb]{$v$}}
\put(1051,-511){\makebox(0,0)[lb]{$u$}}
\put(2026,-1636){\makebox(0,0)[lb]{$t$}}
\put(1276,-736){\makebox(0,0)[lb]{$a_l$}}
\put(1276,-1111){\makebox(0,0)[lb]{$b_l$}}
\put(2026,-736){\makebox(0,0)[lb]{$y_k$}}
\put(2026,-1186){\makebox(0,0)[lb]{$x_k$}}
\put(2851,-961){\makebox(0,0)[lb]{$q_k$}}
\put(1651,-3361){\makebox(0,0)[lb]{$b_lx_k$}}
\put(1651,-2986){\makebox(0,0)[lb]{$y_ka_l$}}
\put(2251,-961){\makebox(0,0)[lb]{$z_k$}}
\put(2251,-3211){\makebox(0,0)[lb]{$z_k$}}
\end{picture}

On this picture $p,u$ are germs of boundary curve  $s_i$, $t,v$ of $s_j$, $q$ of $s_k$.  We remember that $\psi_{m(2)}(\g)=\sum_l a_l \otimes b_l \in \m{i}{A_i}^{op} \otimes \m{j}{A_j}$ and $\psi_{m(3)}(\g)=\sum_k y_k \otimes x_k \otimes z_k \in \m{i}{A_i}^{op} \otimes \m{j}{A_j} \otimes \m{k}{A_k}$. One cam read off the equations from the lower part of the picture:
\begin{equation} \label {E:tgh}
\psi_{m(3)}(\g) \psi_{m(2)}(\g)\otimes 1=\psi_{m(3)}(\g)
\end{equation}
{\bf Interpretation of equation \ref {E:tgh}}
 Fix a fat graph $\g \in M_{g,n}^{2,3}$.  We want to define $\psi_m(\g/e)$ for a fat graph $\g/e$ where $e$ is chosen  as above and equation \ref {E:tgh} satisfies. Denote $Vert(\g)$ the set of vertices of $\g$. Then  $Vert(\g/e) \subset Vert(\g)$ . Define a function $m \rightarrow \psi_m(\g/e)$ as a restriction of  $m \rightarrow \psi_m(\g)$ on $Vert(\g/e)$. Then it is easy to prove that $L(e)\psi(\g)=\psi(\g/e)$.

{\bf Equations for  three-tensors}Suppose we are given two fat graphs $\g, \g' \in Ob(\tmm)$ that are equal, except for parts shown in the upper part of the next picture.   It is clear that if the equal parts of graphs satisfy property described in definition \ref{D:nsk}  defining the set $M_{g,.n}^{2,3}$ one can split in two  an appropriate edge of the graph to make them elements of   $M_{g,.n}^{2,3}$.We know what element of the algebra assign to this two-vertex and  it  presence has no effect on the equation we shall describe now. Therefore  we shall neglect this extra vertex  and treat  $\g, \g'$ as they were elements of $M_{g,.n}^{2,3}$

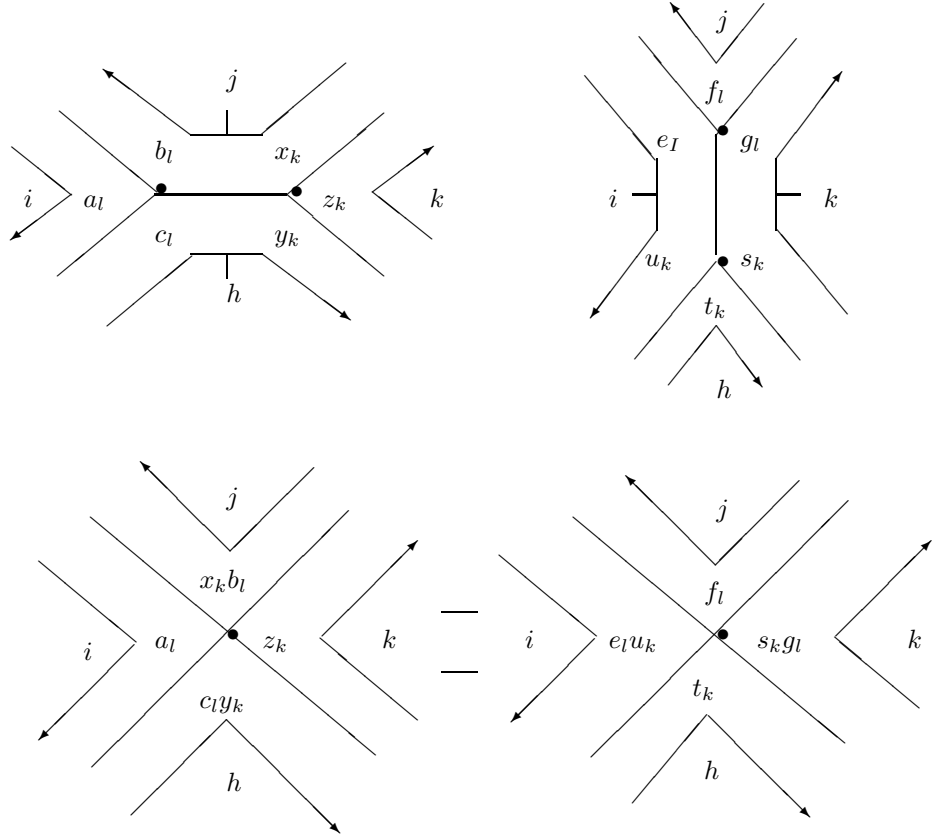
\setlength{\unitlength}{3947sp}%
\begin{picture}(5799,5199)(289,-4423)
\thinlines
\put(1201,-436){$\bullet $}
\put(2048,-458){$\bullet $}
\put(4726,-900){$\bullet $}
\put(4726,-75){$\bullet $}
\put(1651,-3240){$\bullet $}
\put(4726,-3240){$\bullet $}
\put(601, 89){\line( 6,-5){612.295}}
\put(1201,-436){\line(-6,-5){612.295}}
\put(1201,-436){\line( 1, 0){825}}
\put(2026,-436){\line( 6,-5){612.295}}
\put(2026,-436){\line( 6, 5){612.295}}
\put(4726,-61){\line( 0,-1){750}}
\put(301,-136){\line( 5,-4){375}}
\put(676,-436){\vector(-4,-3){384}}
\put(901,-1261){\line( 6, 5){531.148}}
\put(1426,-811){\line( 1, 0){450}}
\put(1876,-811){\vector( 4,-3){552}}
\put(2401,389){\line(-6,-5){531.148}}
\put(1876,-61){\line(-1, 0){450}}
\put(1426,-61){\vector(-4, 3){552}}
\put(5026,764){\line(-4,-5){300}}
\put(4726,389){\vector(-3, 4){288}}
\put(3901,314){\line( 5,-6){442.623}}
\put(4351,-211){\line( 0,-1){450}}
\put(4351,-661){\vector(-3,-4){414}}
\put(5551,-1186){\line(-5, 6){442.623}}
\put(5101,-661){\line( 0, 1){450}}
\put(5101,-211){\vector( 3, 4){414}}
\put(4426,-1636){\line( 4, 5){300}}
\put(4726,-1261){\vector( 3,-4){288}}
\put(4743,-862){\line( 5,-6){510.246}}
\put(4724,-856){\line(-5,-6){510.246}}
\put(5230,547){\line(-5,-6){510.246}}
\put(4752,-58){\line(-5, 6){510.246}}
\put(1651, 89){\line( 0,-1){ 75}}
\put(1651, 14){\line( 0,-1){ 75}}
\put(1651,-61){\line( 0, 1){ 75}}
\put(1651,-811){\line( 0,-1){ 75}}
\put(1651,-886){\line( 0,-1){ 75}}
\put(1651,-961){\line( 0, 1){ 75}}
\put(4351,-436){\line(-1, 0){ 75}}
\put(4276,-436){\line(-1, 0){ 75}}
\put(4201,-436){\line( 1, 0){ 75}}
\put(5251,-436){\line(-1, 0){ 75}}
\put(5176,-436){\line(-1, 0){ 75}}
\put(5101,-436){\line( 1, 0){ 75}}
\put(805,-4036){\line( 1, 1){1612.500}}
\put(794,-2461){\line( 6,-5){1792.623}}
\put(5246,-2236){\line(-1,-1){525}}
\put(4721,-2761){\vector(-1, 1){562.500}}
\put(3826,-2461){\line( 6,-5){1707.934}}
\put(3941,-3974){\line( 1, 1){1612.500}}
\put(1051,-4336){\line( 1, 1){600}}
\put(1651,-3736){\vector( 1,-1){712.500}}
\put(4201,-4261){\line( 5, 6){465.164}}
\put(4676,-3711){\vector( 1,-1){637.500}}
\put(3362,-2702){\line( 6,-5){612.295}}
\put(3962,-3227){\vector(-1,-1){525}}
\put(6001,-3661){\line(-6, 5){531.148}}
\put(5476,-3211){\vector( 1, 1){600}}
\put(2198,-2152){\line(-1,-1){525}}
\put(1673,-2677){\vector(-1, 1){562.500}}
\put(2851,-3736){\line(-6, 5){612.295}}
\put(2251,-3211){\vector( 1, 1){600}}
\put(470,-2740){\line( 6,-5){612.295}}
\put(1070,-3265){\vector(-1,-1){600}}
\put(2941,-718){\line(-5, 4){375}}
\put(2566,-418){\vector( 4, 3){384}}
\put(3001,-3436){\line( 1, 0){225}}
\put(3001,-3061){\line( 1, 0){225}}
\put(1201,-211){\makebox(0,0)[lb]{$b_l$}}
\put(1201,-736){\makebox(0,0)[lb]{$c_l$}}
\put(1951,-211){\makebox(0,0)[lb]{$x_k$}}
\put(1951,-736){\makebox(0,0)[lb]{$y_k$}}
\put(4651,164){\makebox(0,0)[lb]{$f_l$}}
\put(4876,-136){\makebox(0,0)[lb]{$g_l$}}
\put(4651,-1186){\makebox(0,0)[lb]{$t_k$}}
\put(1876,-3286){\makebox(0,0)[lb]{$z_k$}}
\put(4576,-3586){\makebox(0,0)[lb]{$t_k$}}
\put(4651,-2986){\makebox(0,0)[lb]{$f_l$}}
\put(4051,-3286){\makebox(0,0)[lb]{$e_lu_k$}}
\put(4980,-3286){\makebox(0,0)[lb]{$s_kg_l$}}%%%
\put(1480,-3661){\makebox(0,0)[lb]{$c_ly_k$}}
\put(1480,-2911){\makebox(0,0)[lb]{$x_kb_l$}}
\put(1201,-3286){\makebox(0,0)[lb]{$a_l$}}
\put(751,-511){\makebox(0,0)[lb]{$a_l$}}
\put(2251,-511){\makebox(0,0)[lb]{$z_k$}}
\put(4276,-886){\makebox(0,0)[lb]{$u_k$}}
\put(4351,-136){\makebox(0,0)[lb]{$e_I$}}
\put(4876,-886){\makebox(0,0)[lb]{$s_k$}}
\put(376,-511){\makebox(0,0)[lb]{$i$}}
\put(1651,239){\makebox(0,0)[lb]{$j$}}
\put(2926,-511){\makebox(0,0)[lb]{$k$}}
\put(1651,-1111){\makebox(0,0)[lb]{$h$}}
\put(4051,-511){\makebox(0,0)[lb]{$i$}}
\put(3526,-3286){\makebox(0,0)[lb]{$i$}}
\put(751,-3361){\makebox(0,0)[lb]{$i$}}
\put(4726,-1711){\makebox(0,0)[lb]{$h$}}
\put(4651,-4111){\makebox(0,0)[lb]{$h$}}
\put(1651,-4186){\makebox(0,0)[lb]{$h$}}
\put(5401,-511){\makebox(0,0)[lb]{$k$}}
\put(5926,-3286){\makebox(0,0)[lb]{$k$}}
\put(2626,-3286){\makebox(0,0)[lb]{$k$}}
\put(4726,614){\makebox(0,0)[lb]{$j$}}
\put(4726,-2461){\makebox(0,0)[lb]{$j$}}
\put(1651,-2386){\makebox(0,0)[lb]{$j$}}
\end{picture}

We shall  introduce an equation for    $\psi_m(\g)$ and  $\psi_m(\g')$. First of all functions  $m \rightarrow \psi_m(\g)$ and $m \rightarrow \psi_m(\g)$ coincide on the set $Vert(\g) \backslash \{m(3), \bar m(3)\}=Vert(\g') \backslash \{l(3), \bar l(3)\}$.
\\

    By definition 
\begin{align}
&\psi_{m(3)}(\g) = \sum_l a_l \otimes b_l \otimes c_l \in \m{i}{A_i} \otimes \m{j}{A_j}^{op} \otimes \m{h}{A_h}\\
&\psi_{\bar m(3)}(\g)= \sum_k x_k \otimes z_k \otimes y_k \in \m{j}{A_j}^{op} \otimes \m{k}{A_k}^{op} \otimes \m{h}{A_h} 
\end{align}
 Similarly 
\begin{align}
& \psi_{l(3)}(\g')=\sum_l e_l \otimes f_l \otimes g_l \in \m{i}{A_i} \otimes \m{j}{A_j}^{op} \otimes \m{h}{A_k}^{op} \\
&\psi_{\bar l(3)}(\g')=\sum_k u_k \otimes s_k \otimes t_k \in \m{i}{A_i} \otimes \m{k}{A_k}^{op}\otimes \m{h}{A_h} 
\end{align}

Introduce a notation: suppose we are given $n$ associative algebras $B_1, \dots B_n$ and  an element $b=\sum b_s(1) \otimes \dots \otimes b_s(i-1) \otimes b_s(i+1) \dots \otimes b_s(n) \in \bigotimes_{j \neq i} B_j $ . Denote $b^i$ the element $ \sum b_s(1) \otimes \dots \otimes b_s(i-1) \otimes 1 \otimes b_s(i+1) \dots \otimes b_s(n) \in \bigotimes B_j$
Now we are able to introduce third type of equations:
\begin{align}
 \psi_{m(3)}(\g)^3 & \psi_{\bar m(3)}(\g)^1 =\psi_{l(3)}(\g')^4 \psi_{\bar l(3)}(\g')^2 \in  \label {A:dty}\\
& \in  \m{i}{A_i} \otimes \m{j}{A_j}^{op} \otimes \m{k}{A_k}^{op} \otimes \m{h}{A_h} 
\end{align}
One can read off this equation from the lower part of the picture.

{\bf Symmetries} The last condition which we want to impose is nice behavior of $\psi_m(\g)$ with respect to $Aut(\g)$ . Namely we impose :
\begin{equation}
\psi_{T(m)}(\g)=\psi_m(\g), \ T \in Aut(\g)
\end{equation}

\begin{Def} It is easy to see that the collection of function $\psi_m(\g)$ satisfying equations \ref{E:iwf}, \ref{E:tgh}, \ref{A:dty} is closed under direct sum. In other words they form additive category. The Grothendeick group of it is by definition graph K-theory.
\end{Def}

\begin{Prop}
There is a natural way to extend functions $\psi_m{\g}$ on all fat graphs in $\tm$. The collection of elements $\psi(\g)$ functions $\psi_m {\g}$ define satisfy conditions \ref {E:rfy} and  \ref {E:giy}.
\end{Prop}

\subsection{Interpretation of elements of graphical K-groups}
  We start this subsection with the following 
observation which belongs to  mathematical folklore. Suppose that we are given an associative finite
dimensional  algebra $A$. We deform its multiplication law . As a
rule in the process of deformation  algebra becomes ``more
semi-simple''. In particular if the algebra is already semi-simple then
the conventional theory of deformation gives trivial results . If one
still wants to define some sort of ``discrete'' deformation then can
act as follows. Pick an element $$\psi=\sum a_i\otimes
b_i\otimes c_i \in A \otimes A \otimes A$$. Define a new multiplication in $A$ by the rule:
\begin{equation}
x\overset{\psi}{*}y=\sum a_ixb_iyc_i
\end{equation}

Of course for general tensor $\psi$ the new multiplication is no
longer associative. 

Introduce the equations :

1) $\sum a_ja_i \otimes b_i \otimes c_ib_j \otimes c_j=\sum a_i \otimes
b_ia_j \otimes b_j \otimes c_jc_i$

2) $\sum a_ib_i \otimes c_i=\sum a_i \otimes b_ic_i=1 \otimes1$

\begin{Prop}
Suppose that tensor $\psi$ satisfy the condition 1. Then the algebra
$A$
with multiplication $x\overset{\psi}{*}y$ is associative. Denote the
unit of algebra $A$ by $1$.  If $\psi$ is
a solution of the second equation $1$ is the unit for the algebra with
 multiplication $x\overset{\psi}{*}y$.
\end{Prop} 

{\bf Observation}

The elements $a_{i,j}=1\otimes 1$, $b_{1,2,3}=\sum_i a_i\otimes c_i \otimes b_i$ is an element of group
$K^0_{\Bbb{M}_{0,3}}(A)$. Apparently if one impose cyclic symmetry
condition on tensor $\psi$ one can build elements of $K^0$ groups for
all values of $g$ and $n$(the function $\psi_m(\g)$ see equation \ref {E:vjw} is given by the formula $\psi_m(\g)= \sum_i a_i\otimes c_i \otimes b_i$ for a three-vertex $m$ and $\psi_m(\g)=1$ for a two-vertex $m$)  .

We want to propose a method to construct characteristic classes of
deformations in graphical homology. 

The following theorem is true(it is also a part of math folklore):
\begin{Thm}
Suppose we are given an algebra $A$ , $dim A=n^2$. Consider an algebra
$\bold {Mat}_n$. Then there is a element $\psi \in \bold
{Mat}^{\otimes3}$ subject to relations 1,2, such that $A$ is
isomorphic to the matrix algebra with new multiplication
$\overset{\psi}{\times}$. 
\end{Thm}  

We want to use this theorem to build characteristic classes of
deformations. Namely suppose we are given a deformation $A_h$ of algebra $A$
 and suppose that for $h \ne 0$
$A_h\cong\bold {Mat}_n$. Then there is  an element $\psi_h$ which
defines an element of $K$ group and hence characteristic classes.

\subsection{Relation to Galois cohomology}

Suppose that we are given a Galois extension of a field $K \subset
L$ with the Galois group $G$. Denote by $A*$ invertible elements in a ring $A$. There is a complex which computes Galois
cohomology $H^i(G,L^*)$. Here is its  description.
 $C^n\overset{def}{=}(\underbrace{L\underset{K}{\otimes}...\underset{K}{\otimes}L}_{n})^*$. The differential $d=\prod_{i=1}^n\partial_i^{(-1)^i}$.Let $\sum_s a_s(i)\otimes...\otimes a_s(i)\otimes
a_s(i+1)\otimes... \otimes a_s(n))$ be an element of $C^n$.  Define $\partial_i$ by 
\begin{equation}
\begin{split}
& \partial_i(\sum_s a_s(i)\otimes...\otimes a_s(i)\otimes
a_s(i+1)\otimes... \otimes a_s(n))=\\
& = \sum_s a_s(i)\otimes...\otimes
a_s(i) \otimes 1\otimes a_s(i+1)\otimes... \otimes a_s(n)
\end{split}
\end{equation}

It is easy to  see that cocycles in degree 2 map to the
$K^0_{\Bbb{M}_{0,3}}(K)$.Apparently there is a relation of graphical
homology to the Brauer group.

\end{document}